\documentclass[12pt, twoside,a4paper,reqno]{amsart}
\usepackage[a4paper,width=160mm,top=25mm,bottom=30mm]{geometry}
\usepackage{amsmath,amsthm,amssymb}
\usepackage{hyperref}
\usepackage{amsfonts}
\usepackage{enumerate}
\usepackage{graphicx}

\usepackage[normalem]{ulem}
\usepackage[usenames]{color}

\newcommand{\N}{\mathbb{N}}
\newcommand{\Z}{\mathbb {Z}}
\newcommand{\R}{\mathbb {R}}
\newcommand{\C}{\mathbb {C}}
\newcommand{\E}{\mathbb {E}}
\newcommand{\PP}{\mathbb {P}}

\newcommand{\Mcal}{\mathcal{M}}
\newcommand{\Ac}{\mathcal{A}}
\newtheorem{theorem}{Theorem}

\newtheorem{corollary}[theorem]{Corollary}

\newtheorem{lemma}[theorem]{Lemma}

\newtheorem{remark}[theorem]{Remark}

\begin{document}
\title{Transversal fluctuations for a first passage percolation model}
\author{Yuri Bakhtin and Wei Wu}
\address[Yuri Bakhtin]{Courant Institute of Mathematical Sciences, New York University,
251 Mercer st, New York, NY 10012, USA }
\address[Wei Wu]{Courant Institute of Mathematical Sciences, New York University,
251 Mercer st, New York, NY 10012, USA \newline
\& NYU-ECNU Institute of Mathematical Sciences at NYU Shanghai, 3663 Zhongshan Road North, Shanghai 200062,
China.}

\begin{abstract}
We introduce a new first passage percolation model in a Poissonian
environment on $\mathbb{R}^{2}$. In this model, the action of a path depends on
the geometry of the path and the travel time. We prove that the transversal fluctuation exponent 
for point-to-line action minimizers is at least $3/5$. 
\end{abstract}

\maketitle

\section{Introduction}

In this paper, we introduce a new first passage percolation (FPP) model on 
$\mathbb{R}^{2}$. The path action functional that defines the model depends not only on the geometry of the path as a planar set
but also on the speed of travel or, equivalently, the travel time. 

Our model has several technical advantages over the standard FPP model on $\mathbb{Z}^{2}$. As a result, for point-to-line action minimizers, 
we are able to obtain a lower bound for the transversal fluctuation exponent $\xi \geq 3/5$,  an estimate unknown for the standard FPP.

The standard FPP was originally formulated by
Hammersley and Welsh \cite{HW} as a model of the fluid flow in porus medium.
The model is defined on the square lattice~$\mathbb{Z}^{d}$, via a family 
$\left(\tau_e \right)$ of i.i.d.\
random variables with common distribution $F$. For an edge~$e$ of the lattice~$\mathbb{Z}^{d}$, the random variable~$\tau_e$ represents the passage
time through $e$. The random passage time, or distance $T\left( x,y\right) $ between two
vertices $x,y$ is the infimum of $\sum_{e\in r}\tau \left( e\right) $ over
all paths~$r$ connecting~$x$ to~$y$. If $F$ is continuous, then there is
almost surely a unique path $r^{\ast }$ connecting $x$ to $y$ such that $%
T\left( x,y\right) =\sum_{e\in r^{\ast }}\tau \left( e\right) $. Such optimal paths
are called geodesics. The most interesting questions concerning models of this kind can be formulated in terms of asymptotic properties of geodesics between points $x$ and $y$ in the limit $|x-y|\to\infty$.
We refer to a recent survey \cite{survey-50} for the history of the subject.

Many technical difficulties for the standard FPP on $\mathbb{Z}^{d}$ arise due to 
lattice effects or the lack of isotropy. One of them is the limited knowledge about the limiting shape
$B=\lim_{r\to\infty}\frac{B_r}{r}$, where $B_r=\{x\in\mathbb{Z}^{d}: T(0,x)< r\}$.  The existence of this limit, appropriately understood, 
is implied by the shape theorem (\cite{Ric}, see also \cite{Kes}). For a nice
distribution $F$, it is natural to conjecture that the boundary of $B$ is smooth and uniformly curved (as defined in \cite{New}), but this has
not been proved for \textit{any} $F$. Moreover, many important results for
FPP can only be proved based on such uniform curvature assumption. This
includes the bounds on the two exponents $\chi =\chi \left( d\right) $ and $%
\xi =\xi \left( d\right) $. Roughly speaking, $\chi $ and $\xi $ are defined
so that the standard deviation of $T\left( x,y\right) $ is of order $%
\left\vert x-y\right\vert ^{\chi }$, and the fluctuation of the geodesic $%
r^{\ast }$ about the straight line connecting $x$ to $y$ is of order $%
\left\vert x-y\right\vert ^{\xi }$. The standard FPP is believed to belong
to the KPZ universality class \cite{KPZ}, thus in $d=2$ one believes that $%
\chi =1/3$ and $\xi =2/3$. Under certain technical assumptions, it was
proved by Chatterjee that $\chi =2\xi -1$ \cite{Cha}. However, the known
bounds for $\chi $ and $\xi $ are far from optimal. Indeed, the best bounds
for $\sqrt{\text{Var }T\left( x,y\right) }$ are that it has a $\sqrt{\log
\left\vert x-y\right\vert }$ lower bound (by Newman and Piza \cite{NP}), and
sublinear upper bound (by Benjamini, Kalai and Schramm \cite{BKS}, see also 
\cite{DHS} for an extension). With curvature assumptions, it is known that $%
\chi \geq 1/8$ and $\xi \leq 3/4$ (\cite{NP}). For the lower bound for $\xi $%
, the best known result is by Licea, Newman and Piza \cite{LNP}, who defined
another exponent $\xi ^{\prime }$ that should be closely related to $\xi $ and proved that $\xi ^{\prime }\geq 3/5$.
We give a brief discussion  of weaknesses of the exponent $\xi'$ after the statement of our main result, 
Theorem \ref{main}. Here we only note that, due to lattice effects, it is difficult to
adapt the argument in \cite{LNP} for $\xi $ itself, and no such adaptation is known to us.

In this paper, we study a particular FPP model on $\mathbb{R}^{2}$ based on
homogeneous Poisson point process. For every two points in~$\mathbb{R}^{2}$, 
we introduce an analogue of the passage time between them that,
besides the endpoint locations, depends also on one extra parameter, the travel time. 
It is the infimum of a certain action  functional given
by the difference between a quadratic kinetic energy term and the number of
Poisson points touched by the path (see (\ref{A2}) and (\ref{A1}) below for
precise definitions). Bakhtin, Cator and Khanin \cite{BCK} studied a similar
last passage percolation model in $1+1$ dimension, in the context of the 
the Burgers equation with Poissonian forcing. We remark here that
different models of FPP on $\mathbb{R}^{d}$ have been studied by Howard and
Newman (\cite{HN1}, also the survey \cite{HN2}), and by Vahidi-Asl and
Wierman \cite{VW}. Similarly to the models studied in \cite{HN1} and \cite{VW}, our model
has full rotational invariance. 

Our main result states that for point-to-line action minimizers, the transversal
fluctuation exponent $\xi$ satisfies~$\xi\ge 3/5$. This is the first rigorous lower bound
for $\xi $ obtained for any FPP model. We believe that all the existing results
on the exponents $\chi $ and $\xi $ can be adapted to our case (such as
the relation $\chi =2\xi -1$, $\chi \geq 1/8$ and $\xi \leq 3/4$), but we
will not study them in the present paper.

Another motivation for us is the conjecture that for $d=2$, there exist no
doubly infinite geodesics in the standard FPP on $\mathbb{Z}^{d}$\ (a doubly
infinite geodesic is a doubly infinite path such that every finite segment
of the path is a finite geodesic between the end points). It is equivalent
to the statement that the disordered ferromagnetic Ising model on $\mathbb{Z}%
^{2}$ has only two ground states (namely, all $+$ and all $-$, see \cite{LN}%
). This conjecture is partially confirmed in \cite{LN}, where it is shown
(under a curvature assumption on the limit shape) that for Lebesgue almost 
every $\hat{x},\hat{y}\in S^{1}$, there is a.s. no
doubly infinite geodesic that has asymptotic directions $\left( \hat{x},\hat{%
y}\right) $. This result is strengthened by the recent work \cite{DH} that
rules out the existence of the doubly infinite geodesic with \textit{any}
asymptotic directions (assuming that the limit shape boundary is differentiable). 
The conjecture is still open, since there may
exist doubly infinite geodesic with infinite winding number. 
Nonrigourous arguments in \cite{KS} suggest that no doubly infinite
geodesic can exist as long as $\xi >1/2$. Although our bound ($\xi \geq
3/5$) for the FPP model that we study seems to be sufficient, it is still an
interesting open problem to make the heuristic arguments in~\cite{KS} rigorous.

We prove the lower bound for $\xi $ by adapting the martingale inequality
argument in \cite{NP} and \cite{LNP} to the continuous setting. Similar
ideas were also pursued in \cite{Wut} to study the Brownian motion in a
truncated Poissonian potential. One crucial prerequisite to apply this
argument is the \textit{locality} property for standard FPP. That is, if one
slightly perturbs the label $\tau \left( e\right) $ on one edge (retaining
the value of all the other labels), the action of the minimizer will not
change much. 

The locality property is obvious for standard FPP, but
not for our model. The role of small perturbations of the environment in our model is
played by insertions of additional Poisson points. In the situtation where the minimizer contains a long segment
connecting distant Poisson points and another Poisson point is inserted near that segment, modifying
the path in order to include this extra point may lead to a significant change of the path action (since it is quadratic in the total length of the
path). However, we are able to show that a version of the locality property still
holds. Similar ideas were pursued in \cite{BCK} to obtain moment bounds for
action. 

This paper is organized as follows. In Section \ref{model}, we define
precisely the new FPP model and state our main result, Theorem \ref{main}.
Section \ref{shape} contains moment bounds for the action and auxiliary greedy lattice animals estimates. The
lower bound for $\xi $ is obtained by studying the variance of the action
difference of two minimizers, and we detail its proof in Sections~\ref{upppf}
and~\ref{lowpf}. Finally, Section~\ref{sec:appendix} contains a proof of a useful alternative
representation of the action.

\bigskip
\paragraph{\textbf{Acknowledgments:}} We thank Charles Newman for helpful
discussions. The research of W.W. was supported in part by U.S. NSF
grants DMS-1007524 and DMS-1507019. The research of Y.B. was supported in part by NSF grant DMS-1460595.

\section{The Model and Main Results\label{model}}
We work on the plane $\R^2$ which we often identify with $\C$ for convenience.

Let $\omega $ be a locally finite point configuration
sampled from a
homogeneous Poisson point process of unit intensity on $\mathbb{R}^{2}$. This means that 
(i) given a Borel set $B\subset \mathbb{R}^{2}$, the number $\omega \left( B\right)$ of 
configuration points  in~$B$ is a
Poisson random variable with mean $\left\vert B\right\vert $ (the Lebesgue
measure of $B$); (ii) for
disjoint bounded Borel sets $A_{1},...,A_{m}$, the random variables $\omega
\left( A_{1}\right),...,\omega(A_{m})$ are independent. As usual, we identify
locally finite point configurations $\omega $ with integer-valued locally
bounded Borel measures with a unit atom at each point of the configuration.

For any point configuration $\omega$ and any $s>0$ we denote by 
$C_{\omega }(\left[ 0,s\right]:\mathbb{R}^{2})$ the set of
$\mathbb{R}^{2}$-valued 
piecewise linear paths defined
on $\left[ 0,s\right] $ that visit any Poisson point in $\omega $ at most
once. For any $\gamma\in C_{\omega }(\left[ 0,s\right] :\mathbb{R}^{2})$, we define
the following action functional
\[
 A^s(\gamma)= \frac{1}{2}\int_{0}^{s}\left\vert \dot{\gamma} \left( u\right) \right\vert ^{2}du-\omega _{pp}\left( \gamma \right),
\]
where $\omega _{pp}\left( \gamma \right) $ is the number of Poisson points
touched by $\gamma $ and $|\cdot|$ denotes the Euclidean norm. The
first term (\textquotedblleft kinetic energy\textquotedblright ) depends
only on the geometry of the path, whereas the second term (\textquotedblleft
potential energy\textquotedblright ) is responsible for the interaction with
the environment given by the Poisson process.

For any $s>0$, we can define action between a point $x\in \mathbb{R}^{2}$, and a set $S\subset \mathbb{R}^{2}$ by
\begin{equation*}
A^{s}\left( x,S\right) =\inf_{\substack{ \gamma \in C_{\omega }\left( \left[
0,s\right] :\mathbb{R}^{2}\right)  \\ \gamma \left( 0\right) =x,\gamma
\left( s\right) \in S}} A^s(\gamma).
\end{equation*}
We also denote $A^{s}\left( S\right) =A^{s}\left( 0,S\right) $.  
One can write the optimization problem in two
steps: first minimize the velocity conditioned on the point configurations,
then minimize over all such points. This gives another equivalent definition
of the action (see Section~\ref{sec:appendix} for a proof):

\begin{lemma}
\label{Action}We have
\begin{equation}
A^{s}(x,S) =\inf_{\substack{ N\ge 0,\ \left( x_{i}\right)
_{i=0}^{N+1},\ x_{i}\neq x_{j}  \\ x_{0}=x,\ x_{N+1}\in S}}\left\{ \frac{%
\left( \sum_{i=0}^{N}\left\vert x_{i+1}-x_{i}\right\vert \right) ^{2}}{2s}%
-N\right\},  \label{A2}
\end{equation}%
where the infimum is taken over the number $N\in \mathbb{N}\cup \left\{ 0\right\} $,
locations $\left( x_{i}\right) _{i=1}^{N}$ of distinct Poisson points, and the terminal point $x_{N+1}$ in the set $S$.
\end{lemma}

This result shows that it is sufficient to work with paths understood as sequences of 
points $(x_0,x_1,\ldots,x_N,x_{N+1})$, where $x_1,\ldots,x_N$ are distinct Poisson points, assigning action
\[
 A^s(x_0,x_1,\ldots,x_N,x_{N+1})=\frac{L^2(x_0,x_1,\ldots,x_N,x_{N+1})}{2s}-N
\]
to such a path. Here $L(x_0,x_1,\ldots,x_N,x_{N+1})=\sum_{i=0}^{N}\left\vert
x_{i+1}-x_{i}\right\vert$.

We will be mostly interested
in the point-to-line action.
Given a unit vector $u\in S^{1}$, we define $\mathcal{L}_{u%
}=\left\{ au:a\in \mathbb{R}\right\} $. Also, for $z=\left( x,y\right)
\in \mathbb{R}^{2}$ (or equivalently, $x+iy\in \mathbb{C}$), let~$\Lambda
_{z}$ denote the line passing through $z$ that is perpendicular to $\mathcal{%
L}_{z/\left\vert z\right\vert }$. The main object in this paper is the action from $0$ to
$\Lambda _{t}=\Lambda _{t+i0}=\Lambda _{\left(t,0\right) }$:
\begin{equation}
A^{s}\left( 0,\Lambda _{t}\right) =\inf_{\substack{ \gamma \in C_{\omega
}\left( \left[ 0,s\right] :\mathbb{R}^{2}\right)  \\ \gamma \left( 0\right)
=x,\gamma \left( s\right) \in \Lambda _{t}}}A^s(\gamma) .  \label{A1}
\end{equation}
By the nature of the Poisson point process and the continuity of the action
with respect to individual particle locations, it follows that the geodesic, i.e., the minimizer
in (\ref{A1}), is a.s.-unique. 
\bigskip


When studying $A^{s}\left( 0,\Lambda _{t}\right) $, we are
interested in the space-time scaling $s=ct$, for a constant $c>0$. The
competition between the kinetic and potential energy terms depends on $c$.
We expect that there exists $c^{\ast }\in \left( 0,\infty \right) $ with the following properties: 
if $c<c^{\ast }$, then the kinetic energy dominates, and the
minimizer has KPZ\ fluctuations; when $c>c^{\ast }$, the
environment contribution dominates, and the minimizer keeps wandering in order to
collect more Poisson points, which may lead to a larger fluctuation
exponent. In this article, we will focus on the former case  and provide a 
lower bound for the transversal fluctuation exponent of the
point-to-line geodesics for the action defined by (\ref{A1}) (or (\ref{A2}%
)).

Given $w>0$ and $u\in S^{1}$, we define the cylinder $\mathcal{C}_{%
u}\left( w\right) $ symmetric about $\mathcal{L}_{u}$ and of
width $w$ as%
\begin{equation*}
\mathcal{C}_{{u}}\left( w\right) =\left\{ z\in \mathbb{R}^{2}:\text{dist}%
\left( \left\{ z\right\},\mathcal{L}_{{u}}\right) \leq w\right\},
\end{equation*}%
where for $A,B\subset \mathbb{R}^{2}$, dist$\left( A,B\right) :=\inf_{x\in
A}\inf_{y\in B}\left\vert x-y\right\vert $. Given $x\in \mathbb{R}^{2}$, $%
S\subset \mathbb{R}^{2}$, and $t>0$, we denote by $\mathcal{M}\left(
x,S,t\right) \in C\left( \left[ 0,t\right] :\mathbb{R}^{2}\right) $
the path providing the minimal action in the definition of $A^{t}\left(x,S\right) $.
We denote $\mathcal{M}\left( S,t\right) =$ $\mathcal{M}\left(
0,S,t\right) $ and, for $s\in \left[ 0,t\right]$,  use $\mathcal{M}\left( s\right)$ for 
$\mathcal{M}\left(S,t\right) \left( s\right) $, when $S$ and $t$ are clearly defined by the
context. We also denote $\mathcal{M}\left( x,y,t\right) =\mathcal{M}\left(
x,\left\{ y\right\},t\right) $.

Let us now define the transversal fluctuation exponent $\xi $ for $\mathcal{M}$
by%
\begin{equation}
\xi =\sup \left\{ \gamma :\limsup_{t\rightarrow \infty }\sup_{u\in
S^{1}}\mathbb{P}\left( \mathcal{M}\left( \Lambda _{t{u}},ct\right)
\subset \mathcal{C}_{{u}}\left( t^{\gamma }\right) \right) <1\right\} .
\label{xi}
\end{equation}%
Equivalently, we can define $\xi$ by 
\begin{equation}
\xi =\sup \left\{ \gamma :\limsup_{t\rightarrow \infty }\mathbb{P}\left( 
\mathcal{M}\left( \Lambda _{t},ct\right) \subset \mathcal{C}_{e_{1}}\left(
t^{\gamma }\right) \right) <1\right\},  \label{xi1}
\end{equation}
where $e_{1}=\left( 1,0\right) $. This is an immediate consequence of the following result
which in turn follows from rotational invariance of the action functional and 
the distribution of Poisson point process.

\begin{lemma}
\label{rot}Let $x\in \mathbb{R}^{2}$, $S\subset \mathbb{R}^{2}$, $\theta \in
\lbrack 0,2\pi )$, and let $R_{\theta }$ be a rotation defined by $R_{\theta
}\left( S\right) =\left\{ e^{i\theta }y:y\in S\right\} .$ Then
$A^{s}\left( x,\Lambda _{t}\right) =A^{s}\left( R_{\theta }(x),\Lambda
_{e^{i\theta }t}\right)$ and  $\mathcal{M}\left( x,S,t\right) =\mathcal{M}%
\left( R_{\theta }(x),R_{\theta }(S),t\right)$ in distribution.
\end{lemma}

We are ready to state our main result giving a lower bound for $\xi $:

\begin{theorem}
\label{main}There is $c^{\ast }\in ( 0,1] $ such that for
all $c<c^{\ast }$, we have $\xi \geq 3/5.$
\end{theorem}

The proof of Theorem \ref{main} will be given in Sections~\ref{upppf} 
and~\ref{lowpf}. A concrete value of $c^{\ast }$ can be obtained by tracking through
the proofs. Although we
further conjecture $c<1$ suffices for Theorem \ref{main} to hold, we do not
try to optimize the value of $c^{\ast }$.

In the case of standard FPP on $\mathbb{Z}^{2}$,
Theorem \ref{main} is established in \cite{LNP}, for another exponent $\xi
^{\prime }$ closely related to $\xi $. In fact, $\xi ^{\prime }$ defined in 
\cite{LNP} has two weaknesses: it only guarantees $\mathcal{M}\left( \Lambda
_{t_{n}u_{n}},ct_{n}\right) $ not being confined in $\mathcal{C}%
_{u_{n}}\left( t_{n}^{\gamma }\right) $, for some sequences $\left(
u_{n},t_{n}\right) $ such that $t_{n}\rightarrow \infty $; also, the
definition of $\xi ^{\prime }$ involves replacing the minimizer by near
minimizers. By working with a continuum model that has rotational
invariance, we are able to overcome these weaknesses.

Throughout the paper, we write $f(n)\preccurlyeq g(n)$, or $%
f(n)=O(g(n))$ if there exists $C<\infty $, such
that for all $n$, $f(n)\leq Cg(n)$. We write $%
f(n)\asymp g(n)$ if $f(n)\preccurlyeq
g(n)$ and $g(n)\preccurlyeq f(n)$.

\section{Lattice Animals and Moment Bounds for the Action}\label{shape}
We begin with some useful terminology and auxiliary results.

A {\it lattice animal} is a connected subset of $\mathbb{Z}^{2}$ containing the origin. The set of lattice
animals of size $n$ is denoted by $A(n)$. Given a function $X:\Z^2\to\R$, 
the weight of a lattice animal $\mathcal{A}\in
A(n)$, is defined by $N_{\Ac}=N_{\Ac}(X)=\sum_{k\in \Ac}X_k $. The weight of the \textit{greedy lattice
animal} of size $n$ is defined as
\begin{equation}
\label{Nn}
N_{n}=\max_{\Ac\in A(n)} N_{\Ac}(X).
\end{equation}
The following tail bound for the weight of greedy lattice animals is a version of a general estimate
established in \cite{CGGK} (see the remark after (2.12) in \cite{CGGK}) specialized to the Poissonian case.

\begin{lemma}\label{lem:greedytail-general}
Let $(X_j)_{j\in\Z^2}$ be a family of i.i.d.\ Poisson random variables with mean $\lambda>0$.
There is $\rho>0$ such that if $y\ge y_0=(e^3\lambda)\vee \rho$, then  
\[
\PP(N_n> yn)\le e^{-yn},\quad n\in\N.
\]
\end{lemma}
\begin{proof}
We recall that
$\E e^{tX_0}= \exp(\lambda (e^t-1))$ for all $t>0$.
Since the number of lattice animals of size $n$ is bounded by $e^{\rho n}$ for some constant $\rho>0$,
we can use Markov's inequality to write
\[
 \PP(N_n> yn)\le e^{\rho n} \PP  ( X_1+\ldots+X_n>yn)\le e^{f(t,y)n},
\]
where $X_1,\ldots,X_n$ are i.i.d.\ $\lambda$-Poisson, and $f(t,y)=\rho+\lambda(e^t-1)-ty$.
The minimum of $f(t,y)$ over all $t$ is attained at $t=\ln(y/\lambda)$ and equals
$g(y)=\rho-\lambda+y -y \ln(y/\lambda)$.
Clearly, $ g(y)\le \rho +y - 3y \le -y,$
and the lemma follows.
\end{proof}

\begin{corollary}
\label{cor:moments}
Under the conditions of Lemma~\ref{lem:greedytail-general}, $\E N_n^k\preccurlyeq n^k$. 
\end{corollary}
\begin{proof} For large $n$,
\[
 \E N_n^k\le (ny_0)^k\PP(N_n\le ny_0) + \sum_{r> ny_0} (r+1)^k e^{-r}, 
\]
and the desired estimate follows.
\end{proof}


To apply the greedy lattice animal estimates to the study of action minimizers, we partition $%
\mathbb{R}^{2}$ into disjoint union of unit squares 
$B_{(i,j) }=[i-\frac{1}{2},i+\frac{1}{2})\times \lbrack j-\frac{1}{2},j+\frac{1}{2})$,
$(i,j) \in \mathbb{Z}^{2}$, 
and fix an arbitrary total ordering of them. We write $X_{(i,j)}=\omega(B_{(i,j)})$ for the number of Poisson points in $B_{(i,j)}$
and use $N_n$ for the greedy lattice animal weight with respect to thus defined random field $(X_{(i,j)})$ of Poisson random variables
with mean $\lambda=1$.

We say that a path $\gamma$
\textit{passes through} $B_{\left(
i,j\right) }$ if $\gamma(s) \in B_{(i,j) }$ for some~$s$. We say that $\gamma$ \textit{touches} $B_{(i,j) }$ if, 
moreover, $\gamma(s) $ is a
Poisson point in $B_{(i,j) }$ for some $s$.
Let $\mathcal{A}=\mathcal{A}\left(S,ct\right) $
denote the set of unit squares that $\mathcal{M}\left(S,ct\right) $ touches. 
Then $\mathcal{A}$ forms a lattice animal of
size $\left\vert \mathcal{A}\right\vert $.

The following moment estimate is the main result of this section.

\begin{lemma}\label{lem:moment-estimate}
\label{moment}Let  $k\in \mathbb{N}$. Then there is $C_k<\infty$ such that for all $c\in(0,1]$ and sufficiently large $t$, the following holds:
if a set $S\subset \R^2$ contains a point $y$ satisfying $|y|=t$, then
\begin{align*}
\mathbb{E}\left\vert A^{ct}(S) \right\vert
^{k}&\leq C_k c^{-k} t^{k},\\
\mathbb{E}N_{\left\vert \mathcal{A}\right\vert }^k&\leq C_kc^{-k/2}t^k,\\
\mathbb{E}L^{2k}&\le
C_{2k} t^{2k},\\
\E \Ac^{2k}&\le C_{2k} t^{2k}.
\end{align*}
where $L$ denotes the total length of $\mathcal{M}(S,ct)$.
\end{lemma}
\begin{corollary}\label{v} Under the conditions of Lemma~\ref{lem:moment-estimate},
let $v_{t}=L/( ct) $ be the speed of the minimizer.
Then $\mathbb{E}v_{t}^{k}$ is bounded for all $k\in \mathbb{N}$.
\end{corollary}

\begin{proof}[Proof of Lemma \protect\ref{moment}]
The number of Poisson points
touched by $\mathcal{M}\left(S,ct\right) $ does not exceed~$N_{\left\vert \mathcal{A}\right\vert }$.
On the other hand, by a simple
geometric consideration, there exists $c_0>0$, such that if 
$\mathcal{M}\left(S,ct\right) $ touches points from~$n$ unit squares, the
kinetic energy is bounded below by%
\begin{equation*}
\frac{\left( \sum \left\vert x_{i+1}-x_{i}\right\vert \right) ^{2}}{2ct}\geq
c_0\frac{n^{2}}{t},
\end{equation*}%
for all $c\in(0,1]$ and for all large $n$. This implies $A^{ct}(S) \geq
c_0\left\vert \mathcal{A}\right\vert ^{2}/t-N_{\left\vert \mathcal{A}%
_{t}\right\vert }$. Also, by taking a straight path that connects $0$ and $y$
with constant speed $1/c$ and does not collect any Poisson points, we can
bound $A^{ct}(S) $ from above by~ $t/(2c)$. Let us use this to derive an
exponential tail estimate for $\left\vert \mathcal{A}\right\vert $.
Let us define $R=R(c)=(cc_0)^{-1/2}\vee (2c_0^{-1}(e^3\vee \rho))$, where~$\rho$ has been introduced in Lemma~\ref{lem:greedytail-general}.

Since $t/(2c) \geq c_0\left\vert \mathcal{A}\right\vert ^{2}/t-N_{\left\vert \mathcal{A}\right\vert }$, we see
that $n>Rt$ implies
\begin{equation*}
\left\{ \left\vert \mathcal{A}\right\vert \geq n\right\} \subset \bigcup
_{m\geq n}\left\{ N_{m}>c_0\frac{m^{2}}{t}-\frac{t}{2c} \right\} \subset
\bigcup _{m\geq n}\left\{ N_{m}>c_0\frac{m^{2}}{2t} \right\}\subset\bigcup _{m\geq n}\left\{ N_{m}>(e^3\vee\rho) m \right\}.
\end{equation*}%
Estimating the probability of the right-hand side by Lemma~\ref{lem:greedytail-general}, we see that there is a  constant~$M$ 
such that for all $c\in(0,1]$,  sufficiently large $t$, and $n>Rt$, 
\begin{equation*}
\mathbb{P}\left( \left\vert \mathcal{A}\right\vert \geq n\right) \leq
M\exp \left( -c_0n^{2}/(2t) \right) .  
\end{equation*}
For any $k\in\N$, this estimate along with Corollary~\ref{cor:moments} implies 
\begin{eqnarray}
\mathbb{E}N_{\left\vert \mathcal{A}\right\vert }^k
&=&\sum_{n\leq Rt}\mathbb{E}N_{n}^k1_{\left\vert \mathcal{A}\right\vert
=n}+\sum_{n>Rt}\mathbb{E}N_{n}^k1_{\left\vert \mathcal{A}\right\vert =n} 
\notag \\
&\leq &\mathbb{E}N^k_{\left[ Rt\right] }+\sum_{n>Rt}\sqrt{\mathbb{E}N_{n}^{2k}%
}\sqrt{\mathbb{P}\left( \left\vert \mathcal{A}\right\vert =n\right) } 
\notag \\
&\leq &C'\left((Rt)^k+\sum_{n>Rt}n^k\exp \left(
-c_0n^{2}/( 4t) \right)\right)  \notag \\
&\leq &C''c^{-k/2}t^k,  \notag
\end{eqnarray}%
for some constants $C',C''<\infty $.
The lemma now follows from this estimate and the following inequalities: $\left\vert A^{ct}(S) \right\vert \leq \max
\left\{ t/(2c),N_{\left\vert \mathcal{A}\right\vert
}\right\},$
$L^{2}/(2ct)\leq N_{\left\vert \mathcal{A}\right\vert }+t/(2c),$ and
$|\Ac|\le 4 \lceil L\rceil$.
\end{proof}

\section{Transversal Fluctuation Upper Bound\label{upppf}}

Let us fix any $\gamma ^{\prime }>\xi $ and denote $\gamma =\frac{\xi +\gamma
^{\prime }}{2}>\xi $. By the definition (\ref{xi1}) and rotational
invariance, we can choose a subsequence $t_{n}\rightarrow \infty$ such
that for any $\theta \in \lbrack 0,2\pi )$,%
\begin{equation}
\mathbb{P}\left( \mathcal{M}\left( \Lambda _{e^{i\theta
}t_{n}},ct_{n}\right) \subset \mathcal{C}_{e^{i\theta }e_{1}}\left(
t_{n}^{\gamma }\right) \right) \rightarrow 1.  \label{cyl}
\end{equation}

\begin{figure}
\includegraphics[width=10cm]{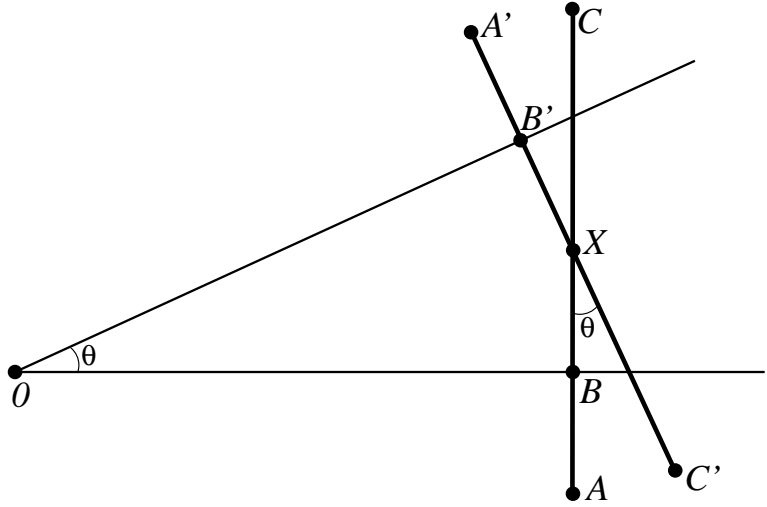} 
\caption{Sets $S(t)$ and $S'(t)$ are represented by line segments $AC$ and $A'C'$. Points $B$ and $B'$ are exactly at distance $t$ from $0$. We have $|AB|=|A'B'|=\frac{1}{2}t^{\gamma'}$ and $|BC|=|B'C'|=\frac{3}{2}t^{\gamma'}$.
We also have asymptotic identities $|BX|=|B'X|\sim \frac{1}{2}t^{\gamma'}$ and $|AX|=|A'X'|\sim |CX|=|C'X|\sim t^{\gamma'}$ as $t\to\infty$. }
\label{fig:S}
\end{figure}

Now, let us define $\theta =\theta_n= t_{n}^{-\left( 1-\gamma ^{\prime }\right) }$ and
introduce the following extensions of segments 
$\mathcal{C}_{e_{1}}\left(t_{n}^{\gamma }\right) \cap \Lambda _{t_{n}}$
and
$\mathcal{C}_{e^{i\theta }e_{1}}\left( t_{n}^{\gamma }\right)
\cap \Lambda _{e^{i\theta }t_{n}}$, see Figure~\ref{fig:S}: 
\begin{eqnarray*}
S\left( t_{n}\right) &=&\left\{ t_{n}e_{1}+ae_{2}:a\in \left[ -\frac{%
t_{n}^{\gamma ^{\prime }}}{2},\frac{3t_{n}^{\gamma ^{\prime }}}{2}\right]
\right\}, \\
S^{\prime }\left( t_{n}\right) &=&\left\{ \left( t_{n}e_{1}+ae_{2}\right)
e^{i\theta }:a\in \left[ -\frac{3t_{n}^{\gamma ^{\prime }}}{2},\frac{%
t_{n}^{\gamma ^{\prime }}}{2}\right] \right\}.
\end{eqnarray*}

The basic idea of our proof is to study the variance of the difference in
action between the point-to-line minimizers connecting the origin to the
lines $\Lambda _{t_{n}}$ and $\Lambda _{e^{i\theta }t_{n}}$. However, we
will study the difference $A^{ct_{n}}\left( S\left( t_{n}\right) \right)
-A^{ct_{n}}\left( S^{\prime }\left( t_{n}\right) \right) $ instead because
it is easier to obtain variance bounds for this quantity.

\begin{remark}\rm
For the rest of the article, we restrict $t$ to the sequence $(t_{n})$ chosen to satisfy~(\ref{cyl}). 
For brevity, we write $t$ for $t_{n}$.
\end{remark}

In the remaining of this section we prove the following proposition.

\begin{lemma} For any $c>0$,
\label{upp}%
\begin{equation*}
\text{Var}\left[ A^{ct}\left( S(t)\right) -A^{ct}\left( S^{\prime }
(t) \right) \right] \preccurlyeq t^{2\left( 2\gamma ^{\prime }-1\right)
}.
\end{equation*}
\end{lemma}

%
%

\begin{proof} The asymptotic identities in the caption of Figure~\ref{fig:S}  imply
that for some constant $c'$ and all sufficiently large $t$, $\sup_{q\in
S( t) }\inf_{y\in S^{\prime }( t) }\left\vert
q-y\right\vert \le c' t^{\gamma ^{\prime }}\theta
\le c' t^{2\gamma ^{\prime }-1}$ and
$\sup_{q\in S^{\prime }( t) }\inf_{y\in S( t)
}\left\vert q-y\right\vert \le c' t^{2\gamma ^{\prime }-1}$.

Given  $\mathcal{M}(S(t),ct) $, we first 
rescale it in time $s\mapsto \frac{ct}{%
ct-c't^{2\gamma ^{\prime }-1}}s$, and then concatenate with a linear path
(with speed bounded by $1$) that connects $\mathcal{M}%
\left( ct\right) $ to $S^{\prime }(t) $ in time bounded by 
$c't^{2\gamma^{\prime }-1}$. The path thus constructed is a suboptimal path for $%
A^{ct}\left( S^{\prime }(t) \right) $, with action at most 
\begin{eqnarray*}
&&A^{ct}(S(t)) +\left( \frac{1}{ct-c't^{2\gamma ^{\prime }-1}}-%
\frac{1}{ct}\right) \frac{L_{t}^{2}}{2}+\frac{c't^{2\gamma^{\prime }-1}}{2} \\
&=&A^{ct}(S(t)) +\frac{c't^{2\gamma ^{\prime }-1}}{\left(
ct-c't^{2\gamma ^{\prime }-1}\right) ct}\frac{L_{t}^{2}}{2}+\frac{c't^{2\gamma^{\prime }-1}}{2}.
\end{eqnarray*}%
This
implies that with probability one, for all $t>0$,
\begin{equation*}
A^{ct}(S(t)) -A^{ct}( S^{\prime }( t))
\geq -\frac{c't^{2\gamma ^{\prime }-1}}{\left( ct-c't^{2\gamma ^{\prime
}-1}\right) ct}\frac{L_{t}^{2}}{2}-\frac{c't^{2\gamma^{\prime }-1}}{2}.
\end{equation*}%
Applying Lemma~\ref{lem:moment-estimate}, we obtain $\mathbb{E}L_{t}^{2}\preccurlyeq
t^{2} $, $\mathbb{E}L_{t}^{4}\preccurlyeq t^{4}$, and therefore,
\begin{equation*}
\mathbb{E}\left( A^{ct}(S(t)) -A^{ct}( S^{\prime }(t) ) \right) _{-}^{2}\preccurlyeq t^{2\left( 2\gamma ^{\prime
}-1\right) },
\end{equation*}%
where $\left( \cdot \right) _{-}$ denotes the negative part. Similarly, one
can prove that  
\begin{equation*}
\mathbb{E}\left( A^{ct}( S(t)) -A^{ct}( S^{\prime }(t) ) \right) _{+}^{2}\preccurlyeq t^{2\left( 2\gamma ^{\prime
}-1\right) },
\end{equation*}
where $\left( \cdot \right) _{+}$ denotes
the positive part, 
and the desired result follows from the fact that $A^{ct}\left( S(t)\right) $
has the same distribution (in particular, the same mean) as $A^{ct}\left(
S^{\prime }(t) \right) $.
\end{proof}

\begin{remark}\rm
The proof of Lemma \ref{upp} uses only rotational 
invariance of the action and the geometry of the objects that we introduced. The definitions of $\xi,\gamma $ or $\gamma
^{\prime }$ were not used.
\end{remark}

\section{Transversal Fluctuation Lower Bound\label{lowpf}}

Let us state the main result of this section.
\begin{theorem}
\label{low}There is $c^{\ast }\in ( 0,1] $, such that for
all $c<c^{\ast }$, and for sufficiently large $t$,%
\begin{equation*}
\text{Var}\left[ A^{ct}(S(t)) -A^{ct}(S^{\prime }(t)) \right] \succcurlyeq t^{1-\gamma }.
\end{equation*}
\end{theorem}

Combining Lemma \ref{upp} and Theorem \ref{low}, we obtain $2\left( 2\gamma
^{\prime }-1\right) \geq 1-\gamma \geq 1-\gamma ^{\prime }$, i.e. $\gamma
^{\prime }\geq \frac{3}{5}$. Since this inequality holds for any $\gamma ^{\prime }>\xi 
$, we conclude that $\xi \geq \frac{3}{5}$, which finishes the proof of Theorem~\ref{main}.

\subsection{Variance decomposition.}
The proof of Theorem \ref{low} is based on a well-known inequality
for functions of independent random variables (cf. eg. \cite{NP} Lemma~2) which follows from
a martingale representation or mutual orthogonality of projections defined by conditional expectations:

\begin{lemma}\label{lem:martingale-inequality}
Let $\left( \Omega,\mathcal{F},\mathbb{P}\right) $ be some probability
space, $T\in L^{2}\left( \Omega,\mathcal{F},\mathbb{P}\right) $ and let $%
\mathcal{G}_{1},\mathcal{G}_{2},...$ be sub-$\sigma $-fields of $\mathcal{F}$%
. If $\left( \mathcal{G}_{i}\right) $ are mutually independent, then 
\begin{equation*}
\text{Var}\left( T\right) \geq \sum_{k}\text{Var}\left( \mathbb{E}\left[ T|%
\mathcal{G}_{k}\right] \right) .
\end{equation*}
\end{lemma}

To apply this inequality, we take $K<\infty $ the value of which will 
be specified later  (see Section \ref{pf1}),
partition $\mathbb{R}^{2}$ into disjoint union of squares $\left\{ B_{\left(
i,j\right) }\right\} _{(i,j) \in \mathbb{Z}^{2}}$ with side
length $K$, such that $B_{(i,j) }$ is centered at $\left(
Ki,Kj\right) $, and fix a total ordering of them. Let $\Omega =\left\{
\omega _{(i,j) }:(i,j) \in \mathbb{Z}^{2}\right\} $
be the space of indexed sequences of Poisson point processes of intensity
one, on the square of side length $K$ (in particular, the number of Poisson points in one square has Poisson distribution with mean $\lambda=K^2$), let $\mathcal{F}$ be the $\sigma -$field
generated by cylindrical sets in $\Omega $, and  let $\mathcal{G}_{\left(
i,j\right) }=\sigma \left( \omega |_{(i,j) }\right) $. 
Then Lemma~\ref{lem:martingale-inequality} yields:
\begin{equation*}
\text{Var}\left[ A^{ct}(S(t)) -A^{ct}(S^{\prime }(t)) \right] \geq \sum_{(i,j) \in \mathbb{Z}^{2}}%
\text{Var}\left( \mathbb{E}\left[ A^{ct}( S(t)) -A^{ct}(S^{\prime }(t)) |\mathcal{G}_{(i,j) }\right]
\right).
\end{equation*}%
For brevity, we will denote $T=A^{ct}( S(t)) $ and $T^{\prime
}=A^{ct}( S^{\prime }(t)) $.

Next, for each $k=(i,j) \in \mathbb{Z}^{2}$ we classify $\omega_{k}=\omega _{(i,j) }$
according to whether $\omega_k$
\textquotedblleft matters" for $T$ and $T^{\prime }$ or not. For each $k\in \mathbb{%
Z}^{2}$, we denote $\hat{\omega}_{k}=\omega |_{\mathbb{R}^{2}\backslash B_{k}}$%
. Given $\hat{\omega}_{k}$, the function that indicates whether the
minimizing path $\mathcal{M}\left( S(t),ct\right) $ touches
some Poisson points in~$B_{k}$ is a function of $\omega _{k}$. Therefore, we
can define%
\begin{eqnarray*}
\mathcal{D}\left( \hat{\omega}_{k}\right) &=&\left\{ \omega _{k}:\text{for }%
\omega =\omega _{k}+\hat{\omega}_{k},\ \mathcal{M}\left( S\left(
t\right),ct\right) \text{ touches }B_{k}\right\}, \\
\mathcal{D}^{\prime }\left( \hat{\omega}_{k}\right) &=&\left\{ \omega _{k}:%
\text{for }\omega =\omega _{k}+\hat{\omega}_{k},\ \mathcal{M}\left(
S^{\prime }(t),ct\right) \text{ touches }B_{k}\right\} .
\end{eqnarray*}%
By the definition of $T$, there is a nonnegative function 
$g$ such that 
\begin{equation*}
T\left( \omega _{k},\hat{\omega}_{k}\right) =T\left( \omega \right) =T\left(
\emptyset,\hat{\omega}_{k}\right) -g\left( \omega _{k},\hat{\omega}%
_{k}\right),
\end{equation*}%
and if $\omega _{k}\notin \mathcal{D}\left( \hat{\omega}_{k}\right) $, then $%
g\left( \omega _{k},\hat{\omega}_{k}\right) =0$. Similarly,
\begin{equation*}
T^{\prime }\left( \omega _{k},\hat{\omega}_{k}\right) =T^{\prime }\left(
\omega \right) =T^{\prime }\left( \emptyset,\hat{\omega}_{k}\right)
-g^{\prime }\left( \omega _{k},\hat{\omega}_{k}\right),
\end{equation*}
for a nonnegative function $g^{\prime }$ satisfying $g^{\prime }\left( \omega _{k},\hat{\omega}_{k}\right) =0$
for $\omega _{k}\notin \mathcal{D}'\left( \hat{\omega}_{k}\right)$.

We can therefore write $\mathbb{E}\left[ T^{\prime }-T|\mathcal{G}_{k}\right]
$ as the following function of $\omega _{k}$ (we denote by $\mathbb{E}_{%
\hat{\omega}_{k}}$ the expectation with respect to the Poisson point process 
$\hat{\omega}_{k}$ in $\mathbb{R}^{2}\backslash B_{k}$):%
\begin{eqnarray*}
&&\mathbb{E}\left[ T^{\prime }\left( \emptyset,\hat{\omega}_{k}\right)
-T\left( \emptyset,\hat{\omega}_{k}\right) |\mathcal{G}_{k}\right] +\mathbb{%
E}_{\hat{\omega}_{k}}\left[ -g^{\prime }\left( \omega _{k},\hat{\omega}%
_{k}\right) 1_{\omega _{k}\in \mathcal{D}^{\prime }\left( \hat{\omega}%
_{k}\right) }+g\left( \omega _{k},\hat{\omega}_{k}\right) 1_{\omega _{k}\in 
\mathcal{D}\left( \hat{\omega}_{k}\right) }\right] \\
&=&\mathbb{E}\left[ T^{\prime }\left( \emptyset,\hat{\omega}_{k}\right)
-T\left( \emptyset,\hat{\omega}_{k}\right) \right] +\mathbb{E}_{\hat{\omega}%
_{k}}\left[ -g^{\prime }\left( \omega _{k},\hat{\omega}_{k}\right) 1_{\omega
_{k}\in \mathcal{D}^{\prime }\left( \hat{\omega}_{k}\right) }+g\left( \omega
_{k},\hat{\omega}_{k}\right) 1_{\omega _{k}\in \mathcal{D}\left( \hat{\omega}%
_{k}\right) }\right] .
\end{eqnarray*}

The first term above is a constant, and we will use the following Lemma 3
from \cite{NP} to remove it.

\begin{lemma}
\label{condition}Suppose $U\in L^{2}\left( \Omega,\mathcal{F},\mathbb{P}%
\right) $, $D_{0}$ and $D_{1}$ are disjoint events in $\mathcal{F}$, then 
\begin{equation*}
\text{Var}\left( U\right) \geq \frac{\mathbb{P}\left( D_{0}\right) \mathbb{P}%
\left( D_{1}\right) }{\mathbb{P}\left( D_{0}\right) +\mathbb{P}\left(
D_{1}\right) }\left( x_{1}-x_{0}\right) ^{2},
\end{equation*}%
where $x_{i}=\mathbb{E}\left( U1_{D_{i}}\right) /\mathbb{P}\left(
D_{i}\right) $ for $i=0,1$.
\end{lemma}

To apply the lemma, we set 
\[
U=\mathbb{E}\left[ T^{\prime }-T|\mathcal{G}_{k}
\right],\quad D_{0}=\left\{ \omega :\omega \left( B_{k}\right) =0\right\},\quad 
D_{1}=\left\{ \omega :2\leq \omega \left( B_{k}\right) \leq b\right\},
\]
where $b<\infty $ is a constant to be chosen later. Noting that
on $D_0$, we have 
$U=\mathbb{E}\left[ T^{\prime }\left( \emptyset,%
\hat{\omega}_{k}\right) -T\left( \emptyset,\hat{\omega}_{k}\right) \right] $%
, we obtain 
\begin{eqnarray}
x_{1}-x_{0} &=&\frac{1}{\mathbb{P}\left( D_{1}\right) }\mathbb{E}_{\omega
_{k}}\mathbb{E}_{\hat{\omega}_{k}}[-g^{\prime }\left( \omega _{k},\hat{\omega%
}_{k}\right) 1_{\omega _{k}\in \mathcal{D}^{\prime }\left( \hat{\omega}%
_{k}\right),2\leq \omega \left( B_{k}\right) \leq b}  \notag \\
&&+g\left( \omega _{k},\hat{\omega}_{k}\right) 1_{\omega _{k}\in \mathcal{D}%
\left( \hat{\omega}_{k}\right),2\leq \omega \left( B_{k}\right) \leq b}].
\label{xdiff}
\end{eqnarray}
To estimate $\left( x_{1}-x_{0}\right) ^{2}$ from below, we need the
following two estimates on $g$ and~$g^{\prime }$.

\begin{lemma}
\label{g'}With probability $1,$ $\left\vert g^{\prime }\left( \omega _{k},%
\hat{\omega}_{k}\right) \right\vert \leq \omega \left( B_{k}\right) $. In particular,
$\left\vert g^{\prime }\left( \omega _{k},\hat{\omega}_{k}\right) \right\vert
\leq b$ almost surely on $\left\{ \omega _{k}\in \mathcal{D}^{\prime }\left( 
\hat{\omega}_{k}\right),2\leq \omega \left( B_{k}\right) \leq b\right\}$.

\end{lemma}

\begin{proof}
By definition, $g^{\prime }\left( \emptyset,\hat{\omega}_{k}\right) =0$.
Now consider adding one point to $B_{k}$. The form of the action (\ref{A2})
suggests, that if the new minimizer collects the new point, the kinetic
energy will increase due to the total length increase (by triangle
inequality). This implies the action $T^{\prime }$ decreases at most by $1$.
The general case follows by induction.
\end{proof}

In the proof of the next lemma, we will switch between two probability
measures. We denote by $\mathbb{P}$ the distribution of Poisson point
process in $\mathbb{R}^{2}$, and $\mathbb{P}_{\omega _{k}}$ the distribution
of Poisson point process in $B_{k}$ (which is also the distribution of a
Poisson point process conditioned on the configuration on 
$\sigma(\omega |_{\mathbb{R}^{2}\setminus B_{k}})$ ).

Before we proceed, let us use Lemma \ref{lem:moment-estimate} to fix $C_1$ such that for every $c\in(0,1]$
and all sufficiently
large $t$, 
\begin{equation}
\mathbb{P}\left( L<C_{1}t\right) >0.99.  \label{L}
\end{equation}

\begin{lemma} \label{g} There exists $c_{1}=c_1(c, K)>0$ 
such that 
\begin{equation*}
\mathbb{E}\left[ g\left( \omega _{k},\hat{\omega}_{k}\right) 1_{\omega
_{k}\in \mathcal{D}\left( \hat{\omega}_{k}\right),2\leq \omega \left(
B_{k}\right) \leq b}\right] \geq \frac{c_{1}}{b}\mathbb{P}\left(
L<C_{1}t,\ \omega _{k}\in \mathcal{D}\left( \hat{\omega}_{k}\right),\ 1\leq
\omega \left( B_{k}\right) \leq b-1\right) .
\end{equation*}
\end{lemma}

\begin{proof}
We will show that one can choose $c_{1}>0$ such that $g\left( \omega _{k},\hat{%
\omega}_{k}\right) \geq 1/2$ on an event that is contained in 
$\left\{ \omega _{k}\in \mathcal{D}\left( \hat{\omega}_{k}\right),\ 2\leq
\omega \left( B_{k}\right) \leq b\right\} $ and has probability at least
$2\frac{c_{1}}{b}\mathbb{P}\left( L<C_{1}t,\ \omega _{k}\in \mathcal{D}\left( 
\hat{\omega}_{k}\right),\ 1\leq \omega \left( B_{k}\right) \leq b-1\right) $.
We denote by $m^{j}$ the Lebesgue measure on $B_{k}^{j}=B_k\times \ldots\times B_k$. Conditioned on any
fixed $\hat{\omega}_{k}$, we can write%
\begin{equation*}
\mathbb{P}_{\omega _{k}}\left( \omega _{k}\in \mathcal{D}\left( \hat{\omega}%
_{k}\right),\ 2\leq \omega \left( B_{k}\right) \leq b\right) =\sum_{j=2}^{b}%
\frac{e^{-K^{2}j}}{j!}m^{j}\left( \left( x_{1},...,x_{j}\right) \in
B_{k}^{j}:\ \sum_{i=1}^{j}\delta _{x_{i}}\in \mathcal{D}\left( \hat{\omega}%
_{k}\right) \right) .
\end{equation*}%
Given $\omega \in \Omega $, we will use the superscript, as in $\mathcal{M}%
^{\left( \omega \right) }$ and $L^{\left( \omega \right) }$, to denote the
dependence of $\mathcal{M}$ and $L$ on $\omega $. Given any $r \in
\left( 0,K\right) $, we define 
\begin{multline*}
E_{j}^{r } =\Bigl\{ \left( x_{1},...,x_{j},y\right) \in
B_{k}^{j+1}:\ \sum_{i=1}^{j}\delta _{x_{i}}\in \mathcal{D}\left( \hat{\omega}%
_{k}\right),\quad \exists h\in \left\{ 1,...,j\right\}  \\
 \text{s.t. }x_{h}\in \mathcal{M}^{\left( \sum_{i=1}^{j}\delta
_{x_{i}}+\hat{\omega}_{k}\right) },\ \left\vert y-x_{h}\right\vert <r
,\ L^{\left( \sum_{i=1}^{j}\delta _{x_{i}}+\hat{\omega}_{k}\right)
}<C_{1}t\Bigr\} \subset B_{k}^{j+1}.
\end{multline*}%
When we add another point $y$ to $B_{k}$, either the new minimizer $\mathcal{%
M}^{\left( \sum_{i=1}^{j}\delta _{x_{i}}+\delta _{y}+\hat{\omega}_{k}\right)
}$ will coincide with $\mathcal{M}^{\left( \sum_{i=1}^{j}\delta _{x_{i}}+%
\hat{\omega}_{k}\right) }$, or it will touch $y$ (in which case it also
touches $B_{k}$). Therefore, if $\left( x_{1},...,x_{j},y\right) \in
E_{j}^{r }$, then $\sum_{i=1}^{j}\delta _{x_{i}}+\delta _{y}\in 
\mathcal{D}\left( \hat{\omega}_{k}\right) $. Hence 
\begin{equation*}
E_{j}^{r }\subset \left\{ \left( x_{1},...,x_{j+1}\right) \in
B_{k}^{j+1}:\ \sum_{i=1}^{j+1}\delta _{x_{i}}\in \mathcal{D}\left( \hat{%
\omega}_{k}\right) \right\} .
\end{equation*}%
Define 
\begin{equation*}
E_{0}^{r }\dot{=}\bigcup_{j=1}^{b-1}\left\{ \omega :\ \omega
_{k}=\sum_{i=1}^{j+1}\delta _{x_{i}},\ \left( x_{1},...,x_{j+1}\right) \in
E_{j}^{r }\right\}.
\end{equation*}%
Since the number of Poisson points in $B_{k}$ is a Poisson random variable
with mean $K^{2}$, 
\begin{eqnarray*}
\mathbb{P}_{\omega _{k}}\left( E_{0}^{r }\right) &=&\sum_{j=2}^{b}\frac{%
e^{-K^{2}j}}{j!}m^{j}\left( E_{j-1}^{r }\right) \\
&\geq &\sum_{j=2}^{b}\frac{e^{-K^{2}j}}{j!}\frac{r^{2}}{4}m^{j}\left(
\left( x_{1},...,x_{j-1}\right) :\ \sum_{i=1}^{j-1}\delta _{x_{i}}\in \mathcal{%
D}\left( \hat{\omega}_{k}\right),\ L^{\left( \sum_{i=1}^{j-1}\delta _{x_{i}}+%
\hat{\omega}_{k}\right) }<C_{1}t\right) \\
&\geq &\frac{1}{e^{K^{2}}b}\sum_{j=1}^{b-1}\frac{e^{-K^{2}j}}{j!}\frac{%
r^{2}}{4}m^{j}\left( \left( x_{1},...,x_{j}\right)
:\ \sum_{i=1}^{j}\delta _{x_{i}}\in \mathcal{D}\left( \hat{\omega}_{k}\right),\ 
L^{\left( \sum_{i=1}^{j}\delta _{x_{i}}+\hat{\omega}_{k}\right)
}<C_{1}t\right) \\
&\geq &\frac{r ^{2}}{4e^{K^{2}}b}\mathbb{P}^{\omega _{k}}\left( \omega
_{k}\in \mathcal{D}\left( \hat{\omega}_{k}\right),\ 1\leq \omega \left(
B_{k}\right) \leq b-1,\ L^{\left( \omega _{k}+\hat{\omega}_{k}\right)
}<C_{1}t\right),
\end{eqnarray*}%
where we used the definition of $E_{j-1}^{r }$ to obtain the first
inequality. Since this inequality holds for \textit{all} $\hat{\omega}_{k}$,
integrating over $\hat{\omega}_{k}$ leads to 
\begin{equation*}
\mathbb{P}\left( E_{0}^{r }\right) \geq \frac{r^{2}}{4e^{K^{2}}b}%
\mathbb{P}\left( \omega _{k}\in 
\mathcal{D}\left( \hat{\omega}_{k}\right),\ 1\leq \omega \left( B_{k}\right) \leq b-1,\ L<C_{1}t\right) .
\end{equation*}%
Now, we show that for $r$ small enough, $g\left( \omega _{k},\hat{%
\omega}_{k}\right) \geq 1/2$ on $E_{0}^{r }$. This will finish the
proof, with $c_{1}=r^{2}/( 8e^{K^{2}}) $.

If $\omega \in E_{0}^{r}$, then there is $j\in \left\{1,...,b-1\right\} $ such that the associated Poisson points satisfy
$\left( x_{1},...,x_{j},y\right) \in E_{j}^{r} $. We claim that if~$r $ is small, then letting the path pick $y$
(immediately after $x_{h}$ in the definition of $E_{j}^{r}$) will
decrease the action at least by $1/2$. To see this, we note that, by the triangle inequality,
the kinetic energy will increase at most by%
\begin{eqnarray*}
\frac{(L+2r )^{2}}{ct}-\frac{L^{2}}{ct} =\frac{4r L}{ ct}
+\frac{4r^{2}}{ct}\leq \frac{4r C_{1}}{c}+\frac{4r^2}{ct} \le \frac{1}{2},
\end{eqnarray*}%
if $r $ is small enough. Since the number of
Poisson points touched is increased by $1$, we have $g\left( \omega _{k},%
\hat{\omega}_{k}\right) \geq 1/2$, and the proof is complete.
\end{proof}

Using (\ref{xdiff}), combining Lemmas \ref{g'} and \ref{g}, denoting $c_{2}=%
\mathbb{P}\left( D_{1}\right) >0$, which depends on $K$ (the size of the
boxes), we have%
\begin{eqnarray*}
x_{1}-x_{0} &\geq &\frac{c_{1}}{bc_{2}}\mathbb{P}\left( \omega _{k}\in 
\mathcal{D}\left( \hat{\omega}_{k}\right),\ 1\leq \omega \left( B_{k}\right)
\leq b-1,\ L<C_{1}t\right) \\
&&-\frac{b}{c_{2}}\mathbb{P}\left( \omega _{k}\in \mathcal{D}^{\prime
}\left( \hat{\omega}_{k}\right),\ 2\leq \omega \left( B_{k}\right) \leq
b\right),
\end{eqnarray*}%
or%
\begin{eqnarray*}
\left\vert x_{1}-x_{0}\right\vert &\geq &\frac{1}{c_{2}}\Big(\frac{c_{1}}{b}%
\mathbb{P}\left( \omega _{k}\in \mathcal{D}\left( \hat{\omega}_{k}\right),\ 
1\leq \omega \left( B_{k}\right) \leq b-1,\ L<C_{1}t\right) \notag \\
&&-b\mathbb{P}\left( \omega _{k}\in \mathcal{D}^{\prime }\left( \hat{\omega}%
_{k}\right),\ 2\leq \omega \left( B_{k}\right) \leq b\right) \Big)_{+}.
\end{eqnarray*}
We can now apply this estimate along with Lemma \ref{condition} and, denoting the area of $\mathcal{C}_{e_{1}}(t^{\gamma })$ by
$\left\vert \mathcal{C}_{e_{1}}(t^{\gamma })\right\vert $, using the monotonicity of $(\cdot)_+$ and Cauchy--Schwarz inequality, obtain
\begin{eqnarray}
&&\sum_{k\in \mathbb{Z}^{2}}\text{Var}\left( \mathbb{E}\left[ T-T^{\prime }|%
\mathcal{G}_{k}\right] \right)  \notag \\
&\geq &\sum_{k\in \mathbb{Z}^{2}:B_{k}\subset \mathcal{C}_{e_{1}}(t^{\gamma
})}\text{Var}\left( \mathbb{E}\left[ T-T^{\prime }|\mathcal{G}_{k}\right]
\right)  \notag \\
&\geq &c_{3}\sum_{k\in \mathbb{Z}^{2}:B_{k}\subset \mathcal{C}%
_{e_{1}}(t^{\gamma })}\Big(\frac{c_{1}}{b}\mathbb{P}\left( \omega _{k}\in 
\mathcal{D}\left( \hat{\omega}_{k}\right),\ 1\leq \omega \left( B_{k}\right)
\leq b-1,\ L<C_{1}t\right)  \notag \\
&&-b\mathbb{P}\left( \omega _{k}\in \mathcal{D}^{\prime }\left( \hat{\omega}%
_{k}\right),\ 2\leq \omega \left( B_{k}\right) \leq b\right) \Big)_{+}^{2}  \notag \\
&\geq &c_{3}b^{2}\left\vert \mathcal{C}_{e_{1}}(t^{\gamma })\right\vert
^{-1}\left( \frac{c_{1}}{b^{2}}\sum_{k\in \mathbb{Z}^{2}:B_{k}\subset 
\mathcal{C}_{e_{1}}(t^{\gamma })}\mathbb{P}\left( \omega _{k}\in \mathcal{D}%
\left( \hat{\omega}_{k}\right),\ 1\leq \omega \left( B_{k}\right) \leq
b-1,\ L<C_{1}t\right) \right.  \notag \\
&&\left. -\sum_{k\in \mathbb{Z}^{2}:B_{k}\subset \mathcal{C}%
_{e_{1}}(t^{\gamma })}\mathbb{P}\left( \omega _{k}\in \mathcal{D}^{\prime
}\left( \hat{\omega}_{k}\right),\ 2\leq \omega \left( B_{k}\right) \leq
b\right) \right) _{+}^{2},  \label{varsum}
\end{eqnarray}%
for some $c_{3}>0$. Since $\left\vert \mathcal{C}_{e_{1}}(t^{\gamma })\right\vert \asymp
t^{1+\gamma }$, if we can prove that for sufficiently small $c$ there are 
values of $b$ and $K$ guaranteeing
\begin{eqnarray}
\frac{1}{t}\liminf_{t\rightarrow \infty }\sum_{k\in \mathbb{Z}%
^{2}:B_{k}\subset \mathcal{C}_{e_{1}}(t^{\gamma })}\mathbb{P}\left( \omega
_{k}\in \mathcal{D}\left( \hat{\omega}_{k}\right),\ 1\leq \omega \left(
B_{k}\right) \leq b-1,\ L<C_{1}t\right) &>&0,  \label{1} \\
\frac{1}{t}\limsup_{t\rightarrow \infty }\sum_{k\in \mathbb{Z}%
^{2}:B_{k}\subset \mathcal{C}_{e_{1}}(t^{\gamma })}\mathbb{P}\left( \omega
_{k}\in \mathcal{D}^{\prime }\left( \hat{\omega}_{k}\right),\ 2\leq \omega
\left( B_{k}\right) \leq b\right) &=&0,  \label{2}
\end{eqnarray}%
 then (\ref{varsum}) implies
\begin{equation*}
\sum_{k\in \mathbb{Z}^{2}}\text{Var}\left( \mathbb{E}\left[ T-T^{\prime }|%
\mathcal{G}_{k}\right] \right) \succcurlyeq \frac{t^{2}}{t^{1+\gamma }}%
=t^{1-\gamma }
\end{equation*}
for those $c$,
which yields Theorem \ref{low}. We study (\ref{1}) and (\ref{2}) separately.

\subsection{Analysis of~(\ref{1})} \label{pf1}
Let us find sufficient conditions on $c,K$, and $b$ guaranteeing~(\ref{1})
by comparing the sum on its left-hand side with the
weight of a lattice path in a Bernoulli site percolation model on~$\mathbb{Z}%
^{2} $.

 Let $W_{t}$ be the number of squares $(B_{k})
_{k\in \mathbb{Z}^{2}}$ touched by $\mathcal{M}\left( S(t)
,ct\right) $, such that $1\leq \omega \left( B_{k}\right) \leq b-1$. We have 
\begin{equation}
\sum_{k\in \mathbb{Z}^{2}:B_{k}\subset \mathcal{C}_{e_{1}}(t^{\gamma })}%
\mathbb{P}\left( \omega _{k}\in \mathcal{D}\left( \hat{\omega}_{k}\right), 
1\leq \omega \left( B_{k}\right) \leq b-1,L<C_{1}t\right) \geq \mathbb{E}%
W_{t}1_{\mathcal{M}\left( S(t),ct\right) \subset \mathcal{C}%
_{e_{1}}(t^{\gamma })}1_{L<C_{1}t}\text{.}  \label{lbd}
\end{equation}%
Now, consider the following Bernoulli site percolation model on $\mathbb{Z}%
^{2}$. Given $\omega $, for any $x\in \mathbb{Z}^{2}$, set $w_{x}=0$ if the
square $B_{x}$ centered at $Kx$ satisfies $\omega \left( B_{x}\right) \in %
\left[ 1,b\right] $, otherwise set $w_{x}=1$. Clearly, $(w_{x})_{x\in \mathbb{Z}^{2}}$ 
are i.i.d.\ random variables, and by
taking $K$ and $b$ large if necessary (to ensure that for a Poisson random
variable $N$ with mean $K^{2}$, $\mathbb{P}\left( 1\leq N\leq b\right) $ is
close to $1$), the connected subsets of $1$s are typically very small.

When this happens, the following adaptation of a result in \cite{CGGK} states
that,
given any lattice path with length $n$, with high probability it only
contains a small fraction of~$1$s. Recall from (\ref{Nn}) that $N_{n}$ is
defined to be the weight of the greedy lattice animal of size $n$ that
contains the origin.

\begin{lemma}
\label{greedy}Let $w_{v}$, $v\in \mathbb{Z}^{2}$ be i.i.d.\ $\left\{
0,1\right\}$-valued random variables with $\mathbb{P}\left( w_{v}=1\right)
=\varepsilon $. Then there is $\tilde c<\infty $ such that for all~$n$, 
\begin{equation*}
\mathbb{P}\left( N_{n}>\tilde cn\varepsilon ^{1/3}\right) \leq e^{-\left( \log
n\right) ^{2}}\text{.}
\end{equation*}
\end{lemma}

\begin{proof}
Let $\tilde{w}_{v}=\varepsilon ^{-1/3}w_{v}$. Since $\mathbb{E}\left\vert 
\tilde{w}_{v}\right\vert ^{3}<\infty $,  the family $(\tilde{w}_{v}) $
satisfies the main condition~(2.4) in~\cite{CGGK}, and it follows from
Proposition 1 in \cite{CGGK} that 
\begin{equation*}
\mathbb{P}\left( \max_{\mathcal{A}\in A(n)}\sum_{v\in \mathcal{A%
}}\tilde{w}_{v}>\tilde cn\right) \leq e^{-\left( \log n\right) ^{2}}
\end{equation*}
for some $\tilde c<\infty$ and all $n$.
\end{proof}

The following is an immediate corollary of Lemma \ref{greedy}.

\begin{corollary}
\label{cor}Let $\tilde c$ be given by Lemma \ref{greedy}. If $\delta >0$ and $K>0$ satisfy 
\begin{equation}
\label{eq:condition-1-on-delta-and-K}
\delta
>\tilde ce^{-K^{2}/3}, 
\end{equation}
then one can choose $b$
so that for all~$n$, $\mathbb{P(}N_{n}>\delta n)\leq e^{-\left( \log
n\right) ^{2}}$.
\end{corollary}

\begin{proof}
This claim follows from Lemma \ref{greedy} and 
\begin{equation*}
\delta \geq \tilde c\left( \mathbb{P}\left( w_{v}=1\right) \right) ^{1/3}=\tilde c\left(
e^{-K^{2}}+\sum_{j>b}\frac{e^{-K^{2}}K^{2j}}{j!}\right) ^{1/3},
\end{equation*}%
which holds true for $\delta >\tilde ce^{-K^{2}/3}$, and all $b$ large enough.
\end{proof}

To apply Corollary \ref{cor}, we need a few more definitions. Given
the minimizer $\mathcal{M}( S(t),ct) $, the lattice
animal \textit{associated with} a path $\mathcal{M}$ is defined to consist
of all $(i,j) \in \mathbb{Z}^{2}$ such that $\mathcal{M}$ 
\textit{touches} $B_{(i,j)}$.
An ordered sequence of (not necessarily distinct) vertices $\left(
y_{1},...,y_{J}\right) \subset \mathbb{Z}^{2}$ is said to be the \textit{%
lattice path traced by }$\mathcal{M}(S(t),ct) $, if 
$\mathcal{M}( S(t),ct) $ sequentially (in time) 
\textit{passes through} $B_{y_{1}},...B_{y_{J}}$. The lattice path
traced by the minimizer always forms a connected subset of $\mathbb{Z}^{2}$,
but it may have self-intersections.

Note that $\mathcal{M}\left( \Lambda _{te_{1}},ct\right) \subset \mathcal{C}%
_{e_{1}}\left( t^{\gamma }\right) $ implies $\mathcal{M}(S(t),ct) 
=\mathcal{M}\left( \Lambda _{te_{1}},ct\right) $, thus $%
\mathcal{M}\left( S(t),ct\right) \subset \mathcal{C}%
_{e_{1}}\left( t^{\gamma }\right) $. Therefore,%
\begin{equation}
\mathbb{P}\left( \mathcal{M}\left( S(t),ct\right) \subset 
\mathcal{C}_{e_{1}}(t^{\gamma })\right) \rightarrow 1\text{ as }t\rightarrow
\infty \text{.}  \label{min}
\end{equation}%
Let us denote by $J$ the size of the lattice path traced by $\mathcal{M}\left(
S(t),ct\right) $. We first notice that $J\geq t/K$. For $n\geq
t/K$, consider the event $H_{n}=\left\{ J=n\right\} $. 

Suppose that $\delta >0$ and $K>0$ satisfy~\eqref{eq:condition-1-on-delta-and-K},
and $b\in\N$ is chosen according to Corollary~\ref{cor}. Then on a good event $G_{n}=H_{n}\cap \left\{
N_{n}\leq \delta n\right\} $ with $\mathbb{P}\left( H_{n}\setminus
G_{n}\right) \leq e^{-\left( \log n\right) ^{2}}$, $\mathcal{M}\left(
S(t),ct\right) $ passes through at least $\left( 1-\delta
\right) n$ squares $\left\{ B_{x}\right\} _{x\in \mathbb{Z}^{2}}$ such that $%
\omega \left( B_{x}\right) \in \left[ 1,b\right] $. Therefore, if $\mathcal{M%
}\left( S(t),ct\right) $ collects Poisson points from at least
a positive fraction of these squares, we will obtain a lower bound 
\begin{eqnarray}
&&\mathbb{E}W_{t}1_{\mathcal{M}\left( S(t),ct\right) \subset 
\mathcal{C}_{e_{1}}(t^{\gamma })}1_{L<C_{1}t}  \notag \\
&\geq &\sum_{n\geq t/K}\mathbb{E}W_{t}1_{\mathcal{M}\left( S( t),ct\right) \subset \mathcal{C}_{e_{1}}(t^{\gamma })}1_{L<C_{1}t}1_{G_{n}} 
\notag \\
&\geq &c_0\sum_{n\geq t/K}n\mathbb{P}\Big(G_{n}\cap \left\{ \mathcal{M}%
\left( S( t),ct\right) \subset \mathcal{C}_{e_{1}}(t^{\gamma
})\right\} \cap \left\{ L<C_{1}t\right\} \Big)  \notag \\
&\geq &c_0\frac{t}{K}\Big[\mathbb{P}\Big(\left( \cup _{n\geq
t/K}H_{n}\right) \cap \left\{ \mathcal{M}\left( S( t),ct\right)
\subset \mathcal{C}_{e_{1}}(t^{\gamma })\right\} \cap \left\{
L<C_{1}t\right\} \Big)-\mathbb{P}\left( \cup _{n\geq t/K}\left(
H_{n}\backslash G_{n}\right) \right) \Big]  \notag \\
&\geq &\frac{c_0}{K}t\left[ 1-\mathbb{P}\Big(\mathcal{M}\left( S(
t),ct\right) \not\subset \mathcal{C}_{e_{1}}(t^{\gamma })\Big)-%
\mathbb{P}\left( L\geq C_{1}t\right) -\sum_{n\geq t/K}e^{-\left( \log
n\right) ^{2}}\right]   \notag \\
&\geq &\frac{c_0}{2K}t,  \label{W}
\end{eqnarray}%
for some $c_0>0$ and all large $t$, where we used $\mathbb{P}\left( \cup
_{n\geq t/K}H_{n}\right) =1$ to obtain the second last inequality, and apply
(\ref{L}) and (\ref{min}) to obtain the last inequality. This yields (\ref{1}%
) due to (\ref{lbd}). To prove that $\mathcal{M}$ collects Poisson points
from sufficiently many squares, we need a geometrical lemma
stating that presence of a long segment in the minimizing path
implies that a certain region in $\mathbb{R}^{2}$ is free of Poisson points.
We  recall that $C_{1}$ is the constant from (\ref{L}) 
and the event ${G}_{m}$ was introduced between~\eqref{min} and~\eqref{W}.

Given $x,y\in \mathbb{Z}^{2}$, we denote $| x-y|_1 =
|x_{1}-y_{1}|+|x_{2}-y_{2}| $. 
\begin{lemma}
\label{geometry} 
Suppose that $\mathcal{M}\left( S(t),ct\right) $
sequentially touches Poisson points $\left( x_{i}\right) $. Let $\left(
y_{i}\right) _{i=1}^{m}\subset \mathbb{Z}^{2}$ be the lattice path traced by 
$\mathcal{M}$, such that $x_{i}\in B_{y_{n_{i}}}$ for some non-decreasing
sequence $\left( n_{i}\right) $.  Then for all
large $t$, on the event ${G}_{m}\cap \left\{ L<C_{1}t\right\}$, relations $n_{i}<j<n_{i+1}$ and 
$|y_{j}-y_{n_{i}}|_1 \wedge |y_{j}-y_{n_{i+1}}|_1 \geq C$
imply $\omega( B_{y_{j}}) =0$, where 
$C=2(\left\lceil4C_{1}c^{-1}K\right\rceil +1)$.
\end{lemma}

\begin{proof} Suppose that for some $j$ satisfying the conditions of the lemma,
there is a Poisson point in $B_{y_{j}}$. Let us prove that
one can decrease the action by collecting
this point. This will contradict the optimality of $\mathcal{M}\left( S(
t),ct\right) $ and finish the proof.

Since $\mathcal{M}\left( S(t),ct\right) $ passes through $%
B_{y_{j}}$, the distance between any Poisson point $z\in B_{y_{j}}$ and the
line connecting $x_{i}$ and $x_{i+1}$ is at most $\sqrt{2}K$. Let $z^{\ast }$
be the projection of $z$ on this line. Denote $s_{1}=\left\vert
x_{i}-z^{\ast }\right\vert $, $s_{2}=\left\vert x_{i+1}-z^{\ast }\right\vert 
$, $L_{i}=L-\left\vert x_{i+1}-x_{i}\right\vert =L-s_{1}-s_{2}$. If $%
s_{1},s_{2}>4C_{1}c^{-1}K^{2}$, where $C_{1}$ is the constant from (\ref{L}%
), by collecting $z$ between the visits at $x_{i}$ and $x_{i+1}$, the
kinetic energy is increased at most by%
\begin{eqnarray*}
&&\frac{\left( L_{i}+\sqrt{s_{1}^{2}+2K^{2}}+\sqrt{s_{2}^{2}+2K^{2}}\right)
^{2}}{ct}-\frac{L^{2}}{ct} \\
&=&\frac{2L_{i}}{ct}\left( s_{1}\sqrt{1+2K^{2}/s_{1}^{2}}+s_{2}\sqrt{%
1+2K^{2}/s_{2}^{2}}-s_{1}-s_{2}\right) \\
&&+\frac{2s_{1}s_{2}(\sqrt{1+2K^{2}/s_{1}^{2}}\sqrt{1+2K^{2}/s_{2}^{2}}%
-1)+4K^{2}}{ct} \\
&\leq &\frac{2L_{i}K^{2}}{ct}\left( 1/s_{1}+1/s_{2}\right) +\frac{2\left(
s_{1}+s_{2}\right) K^{2}}{ct}\left( 1/s_{1}+1/s_{2}\right) +\frac{2K^{4}}{%
s_{1}s_{2}ct}+O\left( \frac{1}{t}\right) \\
&\leq &\frac{2LK^{2}}{ct}\left( 1/s_{1}+1/s_{2}\right) +O\left( \frac{1}{t}%
\right),
\end{eqnarray*}%
\newline
where we use $\sqrt{1+x}-1\leq x/2$ to obtain the second last inequality,
and $L=L_{i}+s_{1}+s_{2}$ to obtain the last inequality. Thus, on the event $%
\left\{ L\leq C_{1}t\right\} $, inequality $s_{1},s_{2}>4C_{1}c^{-1}K^{2}$ implies 
that the
increment of kinetic energy is less than $1,$ so collecting $z$ decreases
the action. Finally, since the number of squares
in $\mathbb{Z}^{2}$ that a line with length $l$ can intersect is at most $%
2\left( \left\lceil l\right\rceil +1\right) $ (both its projections on the $%
x $ and $y$ axis are covered by at most $\left\lceil l\right\rceil +1$
unit segments with integer endpoints), and each box has size $K$, the lemma follows.
\end{proof}

Now, consider the lattice path $(y_{1},...,y_{J})$ traced by $\mathcal{M}%
\left( S(t),ct\right) $, where $J\geq t/K$. Note that this
lattice path may have self-intersections. Lemma \ref{geometry} implies for
any subpath $\left( y_{j},...,y_{j+2C}\right) $ of $\left(
y_{1},...,y_{J}\right) $, such that $w_{y_{j}}=...=w_{y_{j+2C}}=0$, $%
\mathcal{M}\left( S(t),ct\right) $ will collect points from at
least one of $\left\{ B_{y_{k}}\right\} _{k=j}^{j+2C}$. To see this, if $%
\mathcal{M}\left( S(t),ct\right) $ doesn't collect any points
from $\left\{ B_{y_{k}}\right\} _{k=j}^{j+2C}$, let%
\begin{eqnarray*}
a_{1} &=&\max \left\{ a:a<j\text{, and }\mathcal{M}\left( S(t)
,ct\right) \text{ touches }B_{y_{a}}\right\}, \\
a_{2} &=&\min \left\{ a:a>j+2C\text{, and }\mathcal{M}\left( S\left(
t\right),ct\right) \text{ touches }B_{y_{a}}\right\}.
\end{eqnarray*}%
Then $\mathcal{M}$ is a straight line between some point in $B_{y_{a_{1}}}$
and some point in $B_{y_{a_{2}}}$. Thus, it passes through some $B_{y_{l}}$,
with $l\in \left\{ j,...j+2C\right\} $, such that 
$|y_{l}-y_{a_{1}}|_1 \wedge |y_{l}-y_{a_{2}}|_1 \geq
C$, and $\omega \left( B_{y_{l}}\right) >0$. This contradicts Lemma \ref%
{geometry}. Therefore, $\mathcal{M}\left( S(t),ct\right) $ 
\textit{touches} at least one of $\left\{ B_{y_{k}}\right\} _{k=j}^{j+2C}$.

We now partition the lattice path into disjoint unions of subpaths with
length $2C+1$. There are at least $\frac{J}{2C+1}$ of them. On 
$G_{n}$, the cardinality of $\left\{ j\in \left\{ 1,...,n\right\}
:w_{y_{j}}\neq 0\right\} $ is at most $\delta n$. Therefore, at least 
\begin{equation*}
\left( \frac{1}{2C+1}-\delta \right) n
\end{equation*}%
of these subpaths satisfy $w_{y_{j}}=...=w_{y_{j+2C}}=0$. The lattice path
may visit the same square multiple times. However, when $w_{y}=0$, the
minimizer visits $B_{y}$ to pick up points from it at most $b$ times. Thus on 
${G}_{n}$, we have $W_{t}\geq b^{-1}\left( \frac{1}{2C+1}-\delta
\right) n$, which leads to (\ref{W}) with constant $c_0=b^{-1}\left( \frac{%
1}{2C+1}-\delta \right) $.  For this constant to be positive,   we need to require
so that $\delta <\left( 2C +1\right) ^{-1}$, i.e., 
\begin{equation*}
\delta <\left(4\lceil 4C_{1}c^{-1}K\rceil +5 \right) ^{-1}.
\end{equation*}
This condition is implied by
\begin{equation}
\delta <\left(16C_{1}c^{-1}K +9 \right) ^{-1}.
\label{rel}
\end{equation}
Combining conditions \eqref{eq:condition-1-on-delta-and-K} and \eqref{rel}, we arrive to the main conclusion of this section: (\ref{1}) 
is implied by 
\begin{equation*}
  \tilde ce^{-K^{2}/3} <\left(16C_{1}c^{-1}K +9 \right) ^{-1},
\end{equation*}
or, equivalently, by
\begin{equation}
 \label{eq:condition-on-K-no-1}
  c>\frac{16C_1K}{\tilde c^{-1}e^{K^{2}/3}-9},
\end{equation}
because if this condition holds true, then we can find  $\delta$ satisfying \eqref{eq:condition-1-on-delta-and-K} and \eqref{rel},
and $b$ according to Corollary~\ref{cor}.

\subsection{Analysis of (\protect\ref{2})} In this section, we find a condition guaranteeing~\eqref{2} and finish the proof of Theorem~\ref{low}.

We split the sum in (\ref{2}) in two:
\begin{equation*}
\sum_{k\in \mathbb{Z}^{2}:B_{k}\subset \mathcal{C}_{e_{1}}\left( t^{\gamma
}\right) }\mathbb{P}\left( \omega _{k}\in \mathcal{D}^{\prime }\left( \hat{%
\omega}_{k}\right),2\leq \omega \left( B_{k}\right) \leq b\right)
=I_{1}+I_{2},
\end{equation*}%
where $I_1$ is the sum over $k$ satisfying $B_{k}\subset \mathcal{C}_{e_{1}}\left(
t^{\gamma }\right) \cap \mathcal{C}_{e^{i\theta }e_{1}}\left( t^{\gamma
}\right)$ and $I_2$ is the sum over $k$ satisfying $B_{k}\subset \mathcal{C}_{e_{1}}\left(
t^{\gamma }\right) \backslash \mathcal{C}_{e^{i\theta }e_{1}}\left(
t^{\gamma }\right)$ 

We first bound $I_{2}$. Note that $\mathcal{M}\left( \Lambda _{e^{i\theta
}t},ct\right) \subset \mathcal{C}_{e^{i\theta }e_{1}}\left( t^{\gamma
}\right) $ implies $\mathcal{M}\left( S^{\prime }(t),ct\right)
\subset \mathcal{C}_{e^{i\theta }e_{1}}\left( t^{\gamma }\right) $. Let $%
\left\vert \mathcal{A}_{t}^{\prime }\right\vert $ denote the number of
squares touched by $\mathcal{M}\left( S^{\prime }(t),ct\right) $%
. Then, by the Cauchy--Schwarz inequality,%
\begin{eqnarray*}
I_{2} &\leq &\frac{1}{t}\sum_{k\in \mathbb{Z}^{2}:B_{k}\subset \mathcal{C}%
_{e_{1}}\left( t^{\gamma }\right) \backslash \mathcal{C}_{e^{i\theta
}e_{1}}\left( t^{\gamma }\right) }\mathbb{E}1_{\mathcal{M}\left( S^{\prime
}(t),ct\right) \text{ touches }B_{k}} \\
&\leq &\frac{1}{t}\mathbb{E}\left[ \left( \sum_{k\in \mathbb{Z}^{2}}1_{%
\mathcal{M}\left( S^{\prime }(t),ct\right) \text{ touches }%
B_{k}}\right) 1_{\left\{ \mathcal{M}\left( \Lambda _{e^{i\theta
}t},ct\right) \subset \mathcal{C}_{e^{i\theta }e_{1}}\left( t^{\gamma
}\right) \right\} ^{c}}\right]  \\
&\leq &\left[ \mathbb{E}\left( \frac{\left\vert \mathcal{A}_{t}^{\prime
}\right\vert }{t}\right) ^{2}\left( 1-\mathbb{P}\left( \mathcal{M}\left(
\Lambda _{e^{i\theta }t},ct\right) \subset \mathcal{C}_{e^{i\theta
}e_{1}}\left( t^{\gamma }\right) \right) \right) \right] ^{1/2}.
\end{eqnarray*}%
Note that $\left\vert \mathcal{A}_{t}^{\prime }\right\vert $ is different from the $%
\left\vert \mathcal{A}\right\vert $ defined in the proof of Lemma \ref%
{moment}, because it counts the number of $K\times K$ squares. However,
if we let $K$ to be an odd integer, then every unit square we considered 
is a subset of one of the $K\times K$ squares,
so we can claim that $\left\vert \mathcal{A}_{t}^{\prime }\right\vert \leq $ $\left\vert 
\mathcal{A}\right\vert $ a.s., so Lemma \ref{moment} implies $\mathbb{E}%
\left\vert \mathcal{A}_{t}^{\prime }\right\vert ^{2}\preccurlyeq t^{2}$
for any $c$, $b$ and an odd~$K$.
Since we are restricted to the subsequence $\left\{ t_{n}\right\} $ that
satisfies (\ref{cyl}), the definiton~(\ref{xi1}) implies that the last
display above goes to zero as $t_{n}\rightarrow \infty $. Therefore,%
\begin{equation*}
I_{2}/t\rightarrow 0\text{ as }t=t_{n}\rightarrow \infty \text{.}
\end{equation*}
Now, to estimate $I_{1}$, we use the geometric fact that $\mathcal{C}%
_{e_{1}}\left( t^{\gamma }\right) \cap \mathcal{C}_{e^{i\theta }e_{1}}\left(
t^{\gamma }\right) $ is a parallelogram centered at $0$, with diameter $%
O\left( t^{1-\left( \gamma ^{\prime }-\gamma \right) }\right) $. Given $%
\beta\in(0,1)$, let $D_{t^{\beta }}$ be the ball of radius~$t^{\beta }$
centered at $0$, and $U_{t}$ be the number of squares contained in $%
D_{t^{\beta }}$ that are touched by $\mathcal{M}\left( S^{\prime }\left(
t\right),ct\right) $. Clearly, for $\beta >1-\left( \gamma ^{\prime
}-\gamma \right) $ and large $t$, $\mathcal{C}_{e_{1}}\left( t^{\gamma
}\right) \cap \mathcal{C}_{e^{i\theta }e_{1}}\left( t^{\gamma }\right)
\subset D_{t^{\beta }}$, and we have 
\begin{eqnarray*}
&&\sum_{k\in \mathbb{Z}^{2}:B_{k}\subset \mathcal{C}_{e_{1}}\left( t^{\gamma
}\right) \cap \mathcal{C}_{e^{i\theta }e_{1}}\left( t^{\gamma }\right) }%
\mathbb{P}\left( \omega _{k}\in \mathcal{D}^{\prime }\left( \hat{\omega}%
_{k}\right),2\leq \omega \left( B_{k}\right) \leq b\right) \\
&\leq &\sum_{k\in \mathbb{Z}^{2}:B_{k}\subset \mathcal{C}_{e_{1}}\left(
t^{\gamma }\right) \cap \mathcal{C}_{e^{i\theta }e_{1}}\left( t^{\gamma
}\right) }\mathbb{E}1_{\mathcal{M}\left( S^{\prime }(t)
,ct\right) \text{ touches }B_{k}} \\
&\leq &\mathbb{E}U_{t}.
\end{eqnarray*}%

Let us find a condition on $c$ and $K$ guaranteeing $\E U_t=o(t)$ (which, roughly speaking, would mean that the number of squares visited on any microscopic level
is macroscopically negligible) and hence~\eqref{2}.

Given $\omega $, let $\tau =\tau \left( \omega \right) $ and $x=x\left(
\omega \right) $ be the last exit time (and the last exit location,
respectively) for $\mathcal{M}\left( S^{\prime }(t),ct\right) $
to exit $D_{t^{\beta }}$. Namely, $\tau =\sup \left\{ s>0:\mathcal{M}\left(
s\right) \in D_{t^{\beta }}\right\} $, and $x=\mathcal{M}\left( \tau
\right) $. 
Let us fix $\eta \in \left( \delta,1\right) $ and use an argument similar to the
proof of Lemma \ref{moment} to prove that $\mathbb{E}U_{t}1_{\left\{ \tau \leq
t^{\eta }\right\} }\preccurlyeq t^{\eta }$.

Indeed, on $\left\{ \tau \leq t^{\eta }\right\} $, by taking a straight path
connecting $0$ and $\partial D_{t^{\beta }}$ with constant speed $t^{\beta
}/\tau $, we see that $A^{\tau }( 0,x) \leq t^{2\beta }/(2\tau)$. 
Arguing as in the proof of Lemma \ref{moment}, we see that for all $c\in(0,1]$
and all sufficiently large $n$,
\begin{equation}
\label{eq:U-ge-n}
\left\{ U_{t}\geq n\right\} \subset \bigcup _{m\geq n}\left\{
N_{m}>c_0\frac{m^{2}}{\tau} -\frac{t^{2\beta }}{2\tau} \right\}= 
\bigcup _{m\geq n}\left\{
N_{m}>
m\left( c_0\frac{m}{\tau} -\frac{t^{2\beta }}{2\tau m}
\right)\right\},
\end{equation}%
where $N_{m}$ is the weight for the greedy lattice animal of size $m$ 
associated with a family of i.i.d.\ Poisson random variables with
mean $\lambda=K^{2}$ indexed by~$\Z^2$.

Let us recall the definition 
$y_0=y_0(K)=e^3 K^2\vee \rho$ from Lemma~\ref{lem:greedytail-general}.
  If $\tau \leq t^{\eta }$, $m\ge n\ge 2y_0(K)c_0^{-1}t^{\eta }$, 
and $t$ is sufficiently large, then 
\[
 c_0\frac{m}{\tau}> 2\frac{t^{2\beta }}{2\tau m},
\]
and the factor on the right-hand side of~\eqref{eq:U-ge-n} can be
estimated via
\[
c_0\frac{m}{\tau} -\frac{t^{2\beta }}{2\tau m} >y_0(K).
\]
Along with these estimates, Lemma~\ref{lem:greedytail-general} gives for all large $t$ and $n\ge 2y_0(K)c_0^{-1}t^{\eta } $:

\begin{equation*}
\mathbb{P}\left( U_{t}\geq n,\tau \leq t^{\eta }\right)
\le \sum_{m\ge n} \PP\left(N_m\ge \frac{c_0m^2}{2t^\eta}\right)
\le \sum_{m\ge n} \exp\left(-\frac{c_0m^2}{2t^\eta}\right)
\le \bar C \exp
\left( -c_0 n^{2}/t^{\eta }\right),
\end{equation*}%
for some $\bar C<\infty $. So, for all $c$ and $K$, we obtain $\mathbb{E}U_{t}1_{\left\{ \tau \leq t^{\eta
}\right\} }\preccurlyeq t^{\eta }$. 

It remains to find a condition that would guarantee that 
$\mathbb{E}U_{t}1_{\left\{ \tau >t^{\eta }\right\} }=o(t)$. We will show that under a certain requirement on smallness of $c$, 
the event $\left\{ \tau >t^{\eta }\right\} $ is very unlikely,
by arguing that typically inequality $\tau >t^{\eta }$ contradicts the optimality of the path.

Suppose $\tau >t^{\eta }$. Given $\Mcal=\mathcal{M}\left( S^{\prime }\left(
t\right),ct\right) $, we construct a new path $\mathcal{\tilde{M}}$ as
follows. We replace the initial segment of $\mathcal{M}$ by a
straight path connecting $0$ and $x$ with speed $1$, then follow the rest of
the path $\mathcal{M}$ with time rescaled by $s\mapsto \frac{ct-t^{\beta }}{%
ct-\tau }s$. 
We denote by $L_{t}^{\prime }$ and $n_{t}^{\prime }$ the total length
of $\mathcal{M}$ and the
number of Poisson points it touches on time interval $[\tau
,ct]$.
By (\ref{A2}), the action of $\Mcal$  and $\tilde \Mcal$ can be decomposed as
\begin{align*}
A^{ct}(\Mcal)&=A^{\tau }(\Mcal) +\frac{\left( L_{t}^{\prime }\right) ^{2}}{%
ct-\tau }-n_{t}^{\prime }\\
A^{ct}(\tilde \Mcal)&= \frac{t^{\beta }}{2}+\frac{\left( L_{t}^{\prime }\right) ^{2}}{ct-t^{\beta
}}-n_{t}^{\prime }.
\end{align*}%

Thus, replacing $\mathcal{M}$ with $\mathcal{\tilde{M}}$, 
the action decreases at least by%
\begin{eqnarray}
\notag
&&A^{\tau }(\Mcal) +\frac{\left(
L_{t}^{\prime }\right) ^{2}}{ct-\tau }-\frac{t^{\beta }}{2}-\frac{\left(
L_{t}^{\prime }\right) ^{2}}{ct-t^{\beta }} \\
\notag
&=&\left( L_{t}^{\prime }\right) ^{2}\frac{\tau -t^{\beta }}{\left( ct-\tau
\right) \left( ct-t^{\beta }\right) }+A^{\tau }(\Mcal) -\frac{t^{\beta }}{2} \\
&\geq &\frac{\tau }{2c^{2}}+A^{\tau }(\Mcal) -\frac{t^{\beta }}{2},
\label{eq:action-decrease}
\end{eqnarray}%
for all large $t$,
where we use $L_{t}^{\prime }\geq t-t^{\beta }$ and $\tau >t^{\eta }$ to
obtain the last inequality.

We now claim that if $c$ is small, then, unless the number of Poisson points in $D_{t^{\beta }}$
is extremely large, the right-hand side of~\eqref{eq:action-decrease} is positive. This will contradict the
fact that $\mathcal{M}$ is the minimizer. For $y>0$, let $G_{n,y}=\left\{ N_{n}\leq ny\right\} $. 
Lemma~\ref{lem:greedytail-general} implies that
\begin{equation}
\label{eq:G-tail}
 \mathbb{P}\left( G_{n,y_0}^{c}\right)\le e^{- ny_0},
\end{equation}
where $y_0=e^3K^2\vee \rho$.


 Let us denote $E_{n}=\left\{U_t =n\right\} $. By the same geometric
consideration we used in proving Lemma \ref{moment}, there exists $c_{1}>0$,
such that if the restriction of $\Mcal$ onto $[0,\tau]$ touches points from at least $n$ squares with size $K$, its kinetic
energy is bounded below by $c_{1}n^{2}K^{2}/\tau$. So, on $E_{n}\cap G_{n,y_0}$ we have 
\begin{eqnarray*}
A^{\tau }(\Mcal)  &\geq &\frac{%
c_{1}n^{2}K^{2}}{\tau}-N_{n} \\
&\geq &\frac{c_{1}n^{2}K^{2}}{\tau }-y_0n \\
&\geq &\min_{n\in \mathbb{N}}\left\{ \frac{c_{1}n^{2}K^{2}}{\tau }%
-y_0n\right\}  \\
&\geq &-\frac{y_0^{2}}{4K^{2}}\frac{\tau }{c_{1}}.
\end{eqnarray*}%
Thus, if we require
\begin{equation}
\label{eq:requirement-on-c}
\frac{1}{2c^{2}}-\frac{y_0^{2}}{4K^{2}}\frac{1}{c_{1}}>0,
\end{equation}
then, due to~\eqref{eq:action-decrease}, for $t$ sufficiently
large, on $E_{n}\cap G_{n,y_0}$, the action of $\mathcal{\tilde{M}}$ is less than the action of $%
\mathcal{M}$. We conclude that \eqref{eq:requirement-on-c} implies $E_{n}\cap G_{n,y_0}\cap \left\{ \tau >t^{\eta
}\right\} =\emptyset $, for $t$ sufficiently large.
Therefore, assuming that \eqref{eq:requirement-on-c} holds, we use~\eqref{eq:G-tail} to obtain
\begin{eqnarray*}
\mathbb{E} U_t
1_{\left\{ \tau >t^{\eta }\right\} } &\leq &\sum_{n\geq \log t}n\mathbb{P}%
\left( E_{n}\cap \left\{ \tau >t^{\eta }\right\} \right) +\sum_{n\leq \log
t}n \\
&\leq &\sum_{n\geq \log t}n\mathbb{P}\left( \left( G_{n,y_0}\right)
^{c}\right) +\left( \log t\right) ^{2} \\
&\preccurlyeq &\sum_{n\geq \log t}ne^{-ny_0}+\left( \log t\right)
^{2}\preccurlyeq \left( \log t\right) ^{2}.
\end{eqnarray*}%
Rewriting~\eqref{eq:requirement-on-c} as
\begin{equation}
\label{eq:condition-on-c-and-K-no-2}
c<\frac{K\sqrt{2c_1}}{e^3K^2\vee\rho}, 
\end{equation}
we obtain that \eqref{eq:condition-on-c-and-K-no-2} implies~\eqref{2}.

\medskip 

We are now ready to finish the proof of Theorem~\ref{low}. 
\begin{proof}[Proof of Theorem~\ref{low}]
It remains to check that for every sufficiently small $c$ there is $K$ such that~\eqref{eq:condition-on-K-no-1} and~\eqref{eq:condition-on-c-and-K-no-2} both hold.
This is true since the right-hand side of~\eqref{eq:condition-on-K-no-1} decays to zero as a function of $K$ faster than that of~\eqref{eq:condition-on-c-and-K-no-2}. This finishes the proof
of Theorem~\ref{low} (and Theorem~\ref{main}).
\end{proof}

\section{Proof of Lemma~\ref{Action}}
\label{sec:appendix}
Elementary variational calculus implies that the path minimizing the kinetic action between points 
$x$ and $y$ over time $t$ is given by the motion with constant velocity $(y-x)/t$, and the resulting optimal
action is $\frac{|y-x|^2}{2t}$.
Applying this to the definition of the action~(\ref{A1}), we see that
\begin{eqnarray}
\notag
A^{s}\left( x,S\right)  &=&\inf_{\substack{ \gamma \in C_{\omega }\left( %
\left[ 0,s\right] :\mathbb{R}^{2}\right)  \\ 
\gamma \left( 0\right)
=x,\ \gamma \left( s\right) \in S}}\left\{ \frac{1}{2}\int_{0}^{s}\left\vert 
\dot{\gamma}\left( u\right) \right\vert ^{2}du-\omega _{pp}\left( \gamma
\right) \right\}  \\
&=&\inf_{N\geq 0}
\inf_{\substack{(x_{i})_{i=0}^{N+1},\ x_{i}\neq x_{j}  \\ x_{0}=x,\ x_{N+1}\in S}}
\,\inf_{0=t_0<t_1<\ldots< t_{N+1}=s}
\left\{ \frac{1}{2}\sum_{i=0}^{N}\frac{\left\vert
x_{i+1}-x_{i}\right\vert ^{2}}{t_{i+1}-t_{i}}-N\right\},
\label{eq:expr-for-action}
\end{eqnarray}%
where $x_1,\ldots,x_N$ are distinct Poissonian points. 
Fixing $N$ and those points, 
we first minimize over $\left( t_{i}\right)$ by studying the
Lagrangian%
\begin{equation*}
\frac{1}{2}\sum_{i=0}^{N}\frac{\left\vert x_{i+1}-x_{i}\right\vert ^{2}}{%
t_{i+1}-t_{i}}-N+\lambda \left( \sum_{i=0}^{N}\left( t_{i+1}-t_{i}\right)
-s\right),
\end{equation*}%
the stationary point of which satisfies%
\begin{equation*}
\sqrt{\frac{\lambda }{2}}\left( t_{i+1}-t_{i}\right)
=\left\vert x_{i+1}-x_{i}\right\vert,\quad i=0,...,N.
\end{equation*}
Therefore, the speed 
$v=\sqrt{\frac{\lambda }{2}}=|x_{i+1}-x_{i}|/(t_{i+1}-t_i)$ of the optimal path does not depend on $i$. Since 
$\sum \left( t_{i+1}-t_{i}\right) =s$, we have $v=\sum \left\vert x_{i+1}-x_{i}\right\vert /s$. Using this in~\eqref{eq:expr-for-action}, we finish
the proof of the Lemma.

\bibliographystyle{alpha}
\bibliography{FPP}

\begin{thebibliography}{CGGK93}

\bibitem[AHD15]{survey-50}
Antonio Auffinger, Jack Hanson, and Michael Damron.
\newblock 50 years of first passage percolation.
\newblock {\em arXiv:1511.03262}, 2015.

\bibitem[BCK14]{BCK}
Yuri Bakhtin, Eric Cator, and Konstantin Khanin.
\newblock Space-time stationary solutions for the {B}urgers equation.
\newblock {\em Journal of the American Mathematical Society}, 27(1):193--238,
  2014.

\bibitem[BKS03]{BKS}
Itai Benjamini, Gil Kalai, and Oded Schramm.
\newblock First passage percolation has sublinear distance variance.
\newblock {\em The Annals of Probability}, 31(4):1970--1978, 2003.

\bibitem[CGGK93]{CGGK}
J~Theodore Cox, Alberto Gandolfi, Philip~S Griffin, and Harry Kesten.
\newblock Greedy lattice animals {I}: Upper bounds.
\newblock {\em The Annals of Applied Probability}, pages 1151--1169, 1993.

\bibitem[Cha13]{Cha}
Sourav Chatterjee.
\newblock The universal relation between scaling exponents in first-passage
  percolation.
\newblock {\em Annals of Mathematics}, 177(2):663--697, 2013.

\bibitem[DH15]{DH}
Michael Damron and Jack Hanson.
\newblock Bigeodesics in first-passage percolation.
\newblock {\em arXiv:1512.00804}, 2015.

\bibitem[DHS13]{DHS}
Michael Damron, Jack Hanson, and Philippe Sosoe.
\newblock Sublinear variance in first-passage percolation for general
  distributions.
\newblock {\em Probability Theory and Related Fields}, pages 1--36, 2013.

\bibitem[HN97]{HN1}
C~Douglas Howard and Charles~M Newman.
\newblock Euclidean models of first-passage percolation.
\newblock {\em Probability Theory and Related Fields}, 108(2):153--170, 1997.

\bibitem[HN01]{HN2}
C~Douglas Howard and Charles~M Newman.
\newblock Geodesics and spanning trees for {E}uclidean first-passage
  percolation.
\newblock {\em Annals of Probability}, pages 577--623, 2001.

\bibitem[HW65]{HW}
John~M Hammersley and DJA Welsh.
\newblock First-passage percolation, subadditive processes, stochastic
  networks, and generalized renewal theory.
\newblock In {\em Bernoulli 1713 Bayes 1763 Laplace 1813}, pages 61--110.
  Springer, 1965.

\bibitem[Kes86]{Kes}
Harry Kesten.
\newblock Aspects of first passage percolation.
\newblock In {\em {\'E}cole d'{\'E}t{\'e} de Probabilit{\'e}s de Saint Flour
  XIV-1984}, pages 125--264. Springer, 1986.

\bibitem[KPZ86]{KPZ}
Mehran Kardar, Giorgio Parisi, and Yi-Cheng Zhang.
\newblock Dynamic scaling of growing interfaces.
\newblock {\em Physical Review Letters}, 56(9):889, 1986.

\bibitem[KS91]{KS}
Joachim Krug and Herbert Spohn.
\newblock Kinetic roughening of growing surfaces.
\newblock {\em C. Godreche, Cambridge University Press, Cambridge}, 1(99):1,
  1991.

\bibitem[LN96]{LN}
Cristina Licea and Charles~M Newman.
\newblock Geodesics in two-dimensional first-passage percolation.
\newblock {\em The Annals of Probability}, 24(1):399--410, 1996.

\bibitem[LNP96]{LNP}
Cristina Licea, Charles~M Newman, and Marcelo~ST Piza.
\newblock Superdiffusivity in first-passage percolation.
\newblock {\em Probability Theory and Related Fields}, 106(4):559--591, 1996.

\bibitem[New95]{New}
Charles~M Newman.
\newblock A surface view of first-passage percolation.
\newblock In {\em Proceedings of the International Congress of Mathematicians},
  pages 1017--1023. Springer, 1995.

\bibitem[NP95]{NP}
Charles~M Newman and Marcelo~ST Piza.
\newblock Divergence of shape fluctuations in two dimensions.
\newblock {\em The Annals of Probability}, pages 977--1005, 1995.

\bibitem[Ric73]{Ric}
Daniel Richardson.
\newblock Random growth in a tessellation.
\newblock In {\em Mathematical Proceedings of the Cambridge Philosophical
  Society}, volume~74, pages 515--528. Cambridge Univ Press, 1973.

\bibitem[VAW90]{VW}
Mohammad~Q Vahidi-Asl and John~C Wierman.
\newblock First-passage percolation on the {V}oronoi tessellation and
  {D}elaunay triangulation.
\newblock In {\em Random graphs}, volume~87, pages 341--359, 1990.

\bibitem[Wut98]{Wut}
Mario~V Wuthrich.
\newblock Superdiffusive behavior of two-dimensional {B}rownian motion in a
  {P}oissonian potential.
\newblock {\em Annals of {p}robability}, pages 1000--1015, 1998.

\end{thebibliography}

\end{document}